\documentclass{m2an}
\pdfoutput=1
\usepackage[T1]{fontenc}
\usepackage[utf8]{inputenc}
\usepackage{lmodern}
\usepackage[english]{babel}

\usepackage{moreverb}
\usepackage{multicol, caption}
\usepackage{booktabs}
\usepackage{graphicx}
\usepackage{subfig}
\usepackage{enumerate}
\usepackage[usenames,dvipsnames]{xcolor}
\usepackage{xfrac}
\usepackage{footnote}
\usepackage{textcomp}
\usepackage{dsfont}
\usepackage{tabularx}
\usepackage{slashbox}
\usepackage{vmargin}
\usepackage{rotating}
\usepackage{amssymb}
\usepackage{amsmath}
\usepackage{mathabx}
\usepackage{amsthm}
\usepackage{mathtools}
\usepackage{lscape}
\usepackage{tikz}
\usepackage{appendix}
\usepackage[nosepfour,np]{numprint}
\npthousandsep{\,}\npdecimalsign{.}
\usepackage{mathrsfs}
\usepackage{placeins}

\usepackage{mathtools}
\mathtoolsset{showonlyrefs=true}

\usepackage{graphicx,color,psfrag,tikz,xcolor,pict2e}
\usetikzlibrary{calc,decorations.markings,positioning,shapes,positioning,arrows,patterns}
\tikzstyle{abstractbox} = [draw=black, fill=white, rectangle,inner sep=11pt, style=rounded corners]
\tikzstyle{abstracttitle}=[fill=white]
\newcommand{\mybox}[3][fill=white]{
    \begin{center}
      \begin{tikzpicture}
        \node [abstractbox, #1] (box) {\begin{minipage}{0.9\linewidth}  #2\end{minipage}};
        \node [abstracttitle, right=11pt] at (box.north west) {#3};
      \end{tikzpicture}
      
    \end{center}
}

\newcommand{\figref}[1]{Fig.~\ref{#1}}

\usepackage[pdftex]{hyperref} 
\hypersetup{
	colorlinks,
	citecolor=black,
	filecolor=black,
	linkcolor=black,
	urlcolor=black
}

\newcommand{\Q}{Q^{ad}}
\renewcommand{\u}{\mathbf{u}}
\renewcommand{\v}{\mathbf{v}}
\newcommand{\z}{\mathbf{z}}
\newcommand{\du}{\mathbf{\delta u}}
\newcommand{\tu}{\mathbf{\tau u}}
\newcommand{\dq}{{\delta q}}
\newcommand{\qs}{{q_\sigma}}
\newcommand{\weak}{\rightharpoonup}

\renewcommand{\O}{{\Omega_0}}
\newcommand{\Oq}{{\Omega_q}}
\renewcommand{\hat}[1]{\widehat{#1}}
\renewcommand{\tilde}[1]{\widetilde{#1}}
\renewcommand{\check}[1]{\widecheck{#1}}
\newcommand{\beq}{\begin{equation}}
\newcommand{\eeq}{\end{equation}}

\newcommand{\fenics}{\emph{FEniCS}}
\renewcommand{\div}{\mathrm{div}}
\newcommand{\cof}{\mathrm{cof}}

\graphicspath{{./}}
\DeclareGraphicsExtensions{.pdf,.ps,.eps,.png,.jpeg,.mps}

\theoremstyle{plain}
\newtheorem{ssmptn}[thrm]{Assumption}

\newenvironment{pbdef}{%
\noindent\ignorespaces%
\setlength{\leftskip}{1cm}%
\em
}

\begin{document}
\bibliographystyle{plain}
\title{Shape optimization for Stokes flow:\\ a reference domain approach}
\author{Ivan Fumagalli, Nicola Parolini \and Marco Verani}\address{MOX - Modellistica e Calcolo Scientifico, Dipartimento di Matematica ``F. Brioschi'', Politecnico di Milano, via Bonardi 9, 20133 Milano, Italy, e-mail: \url{ivan1.fumagalli@mail.polimi.it}, \url{nicola.parolini@polimi.it}, \url{marco.verani@polimi.it}}
%
\date{\today}
\begin{abstract} In this paper we analyze a \emph{shape optimization} problem, with Stokes equations as the state problem, defined on a domain with a part of the boundary that is described as the graph of the control function. The state problem formulation is mapped onto a reference domain, which is independent of the control function, and the analysis is mainly led on such domain.
The existence of an optimal control function is proved, and optimality conditions are derived. After the analytical inspection of the problem, finite element discretization is considered for both the control function and the state variables, and a priori convergence error estimates are derived. Numerical experiments assess the validity of the theoretical results. 
\end{abstract}
%
%
%
%
\maketitle
\section*{Introduction}

Optimal control for partial differential equations \cite{Lions1971} is a challenging field of applied mathematics, thanks to its combination of sophisticated theoretical tools and interesting engineering applications.
Among optimal control problems, shape optimization \cite{Delfour2011,MR1179448} has recently undergone a renewal of interest, mainly due to the wide range of industrial and real world applications, like fluid dynamics \cite{Gunzburger2003} and structural mechanics \cite{MR2270119}, and to the increased computational power available for numerical simulations. Shape optimization aims at finding the solution of problems of the following general form: 
\beq
	\min_{\Omega\in \mathcal{O}} J(\Omega,S(\Omega)),\text{ \ subject to a differential problem $L(S(\Omega))=0$ in $\Omega$},
\eeq
where $J$ is a cost functional, defined on a suitable set $\mathcal{O}$ of admissible domains, $L$ is a differential operator and $S$ is the operator mapping an admissible domain $\Omega\in\mathcal{O}$ to the corresponding solution of the differential problem $L(S(\Omega))=0$ in $\Omega$.







This kind of problems has been widely discussed in the literature, employing different techniques in the description of the set $\mathcal{O}$, generally considered as a proper subset of finite (see, e.g., \cite{Belahcene2003,Ballarin2013}) or infinite (see, e.g., \cite{MR1179448,Delfour2011}) dimensional spaces. The present paper belongs to the latter category, as the boundary of the admissible domains (or a subset of it) is described by the graph of a suitable control function. This approach has been widely adopted by many authors (see, e.g., \cite{MR1969772,haslinger1988finite,begis-glowinski, Eppler:2000,Gunzburger2000,laumen, Kinigera, ABV:2013}). 

Concerning the numerical solution of shape optimization problems, a standard technique is represented by gradient type iterative algorithm, in which the state problem is solved on differently shaped domains at each iteration (see, e.g., \cite{MR2270119, Dogan2007}). A critical point of this approach is the repeated deformation of the computational mesh, leading to an increase of the computational effort and to the possible generation of highly skewed mesh elements. In order to avoid such problems, in this work the reference-domain approach introduced in \cite{Kinigera} is followed, mapping the actual domain and the whole optimization problem onto a reference domain $\Omega_0$. Exploiting this mapping, a priori estimates for the discretization error of the optimization problem are derived, and these results are assessed through numerical tests. Discretization of shape optimization problems and convergence issues have been discussed in other works, such as \cite{MR1969772,haslinger1988finite,MR1993939,Chenais2006a}. However, to the best of our knowledge, only \cite{Eppler2007,Kinigera} provide a convergence rate for the discretization error for the Poisson equation. In this paper we obtain similar convergence results for the Stokes problem; this seems to be the first convergence result for shape optimization problems governed by this class of equations.  

The present paper is organized as follows.
In Section 1, we present the shape optimization problem governed by Stokes equations, and we reformulate it on the reference domain. Within this framework, 
the existence of an optimal solution to the minimization problem is proved. Finally, we consider first order optimality conditions and we provide a boundary-integral expression for them.
Section 2 is devoted to the proof of a priori error estimates for the numerical discretization error of the optimization problem. 
Finally, in Section 3 we present some numerical tests, assessing the theoretical results.
In Appendix A, we discuss the regularity assumptions needed by the a priori estimates, whereas in Appendix B some technical results are proved.


\section{The optimal control problem}\label{sec:def}

The aim of the present paper is to study a shape optimization problem governed by Stokes equations, which reads as follows
\beq
\begin{split}
	\min_{q\in\Q} J(q,\u,p) \text{ \ subject to the following generalized Stokes system:}\\
\left\{
	\begin{aligned}
	\eta \mathbf{u} - \div(\nu\nabla\mathbf{u}) + \nabla p &= \mathbf{f}, \qquad &\text{in } \Omega_q,\\
	\div\ \mathbf{u} &= 0, \qquad &\text{in }\Omega_q,\\
	\mathbf{u} &= 0, \qquad &\text{on } \Gamma_q,\\
	\nu\partial_\mathbf{n }\mathbf{u} - p \mathbf{n} &= \mathbf{g}_N,\qquad &\text{on } \Gamma_1,\\
	\partial_\mathbf{n} u_x = 0,\quad u_y&=0,\qquad &\text{on } \Gamma_2,\\
	\mathbf{u} &= \mathbf{g}_D,\qquad &\text{on } \Gamma_3,
	\end{aligned}
\right.
\end{split}
\label{eq:Stokes}
\end{equation}
where $J$ is a given cost functional to be optimized, $\u=(u_x,u_y)$ and $p$ are the so-called state variables and $q$ is the control function (belonging to the admissible set $\Q$) that identifies the domain 
$\Omega_q$.

\begin{figure}[b]
	\centering
	\includegraphics[width=0.7\textwidth]{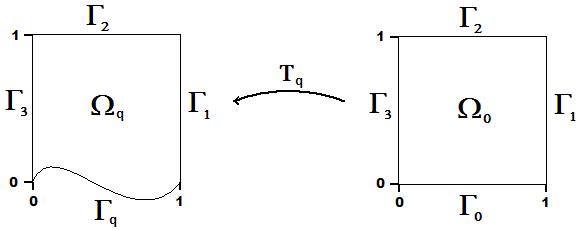}
	\caption{Physical and reference domains}
\label{fig:domini}
\end{figure}

In particular, the control function $q:I=(0,1)\to\mathds{R}$ describes the lower part $\Gamma_q$ of the boundary of domain $\Oq =\{(x,y)\in \mathds{R}^2\ |\ x\in I\ ,\ y\in (q(x),1)\}$. As shown in \figref{fig:domini}(left), the boundary of $\Omega_q$ is partitioned as $\partial\Omega_q=\Gamma_q\cup\Gamma_1\cup\Gamma_2\cup\Gamma_3$.\\
In order to avoid domain degeneration, we fix $\varepsilon\in(0,1)$ a priori, and we introduce the following intermediate set of admissible controls
\beq
	\overline{Q}^{ad}=\{q\in H^3(I)\cap H^1_0(I)\colon q(x)\leq 1-\varepsilon, \ \forall x\in I \}.
\label{eq:Qad1}
\eeq
In the following, it will be useful to have the admissible controls in a bounded set, so we fix a constant $C>0$ and reduce $\overline{Q}^{ad}$ to the following set:
\beq
	Q^{ad} = \{q\in \overline{Q}^{ad}\colon \|q\|_{H^3(I)}\leq C\}.
\label{eq:Qlim}
\eeq
From the above definition, it follows that all the feasible domains $\Omega_q$ are contained in a bounded, convex, hold-all domain $\hat{\Omega}\subset\mathds{R}^2$.

The weak formulation of problem \eqref{eq:Stokes} reads:

\begin{pbdef}
	Find $\mathbf{u} = \hat{\u}+\widetilde{\mathcal{R}}\mathbf{g}_D \ ,\ \hat{\u}\in V_q$ and $p\in P_q$ such that \\
	\beq\left\{
	\begin{aligned}
		a_q(\hat{\u},\mathbf{v}) + b_q(\mathbf{v},p) &= F_q(\mathbf{v}),\qquad&\forall\ \mathbf{v}\in  V_q,\\
		b_q(\hat{\u},\pi) &= -b_q(\widetilde{\mathcal{R}}\mathbf g_D,\pi),
				\qquad&\forall\ \pi \in P_q,\\
	\end{aligned}\right.\\
	\label{eq:Stokesdeb}
	\eeq
	where
	\begin{equation}\begin{split}
		V_q&=\{\v\in[H^1(\Oq)]^2\colon \v=(v_x,v_y)=\mathbf 0\text{ on }\Gamma_3\cup\Gamma_q\text{ and }v_y=0\text{ on }\Gamma_2\},\\
		P_q&=L^2(\Oq),
	\label{eq:spacesq}
	\end{split}\end{equation}
	and
	\beq\begin{split}
	a_q(\mathbf{u},\mathbf{v}) &=\int_{\Omega_q}{\eta \mathbf{u}\mathbf{v}+\nu\nabla \mathbf{u}:\nabla\mathbf{v}},\\
		b_q(\mathbf{v},\pi) &= -\int_{\Omega_q}{\pi\ \div\ \mathbf{v}},\\
		F_q(\mathbf{v}) &= \int_{\Omega_q}{ \mathbf{f}\cdot \mathbf{v}}-a_q(\widetilde{\mathcal{R}}\mathbf{g}_D,\mathbf{v}) + \int_{\Gamma_1}{\mathbf{g}_N\cdot\mathbf{v}\,d\Gamma}.\\
	\end{split}\eeq
\end{pbdef}
Data functions $\eta,\nu,\mathbf{f}$ are defined on the hold-all domain $\hat{\Omega}$,%
\footnote{If not necessary, no special notation will be used to point out whether the entire functions are to be considered, or their restrictions to $\Omega_q$: the distinction will be inferable from the context.}
 boundary data $\mathbf g_N,\mathbf g_D$ are defined on the fixed edges $\Gamma_1,\Gamma_3$, respectively, and $\widetilde{\mathcal{R}}\mathbf{g}_D$ is a continuous lifting of $\mathbf{g}_D$ on $\Omega_q$.
\\

\begin{rmrk}[Well-posedness of the state problem]
Using classical results on Stokes problem (see, e.g., \cite{Girault1986}), we can ensure the well-posedness of \eqref{eq:Stokes}. About data functions, we have to assume what follows:
\footnote{If a particular $q$ is fixed, the conditions need only to be respected on $\Omega_q$. However, in order to be free from dependence on the control, we formulate them on the hold-all domain $\hat{\Omega}$.}
\begin{itemize}
	\item $\nu(\mathbf x)\geq\nu_0>0\ \ \forall \mathbf x\in\hat{\Omega}$,
	\item $\nu,\eta\in L^\infty(\hat{\Omega})$,
	\item $\mathbf{f}\in [H^{-1}(\hat{\Omega})]^2,\ \mathbf g_D\in [H^{1/2}(\Gamma_3)]^2,\ \mathbf g_N\in [H^{-1/2}(\Gamma_1)]^2$.
\end{itemize}
Under these conditions, the following stability estimate holds:
\beq
	\|\tilde{\u}\|_{V_q}+\|\widetilde{p}\|_{P_q} \leq c(\|\mathbf f\|_{[H^{-1}(\hat{\Omega})]^2}+\|\mathbf g_D\|_{[H^{1/2}(\Gamma_3)]^2}+\|\mathbf g_N\|_{[H^{-1/2}(\Gamma_3)]^2}).
	\label{eq:stimauBrezzi}
\eeq
We remark that constant $c$ in \eqref{eq:stimauBrezzi} is independent of $q$, 
since the inf-sup constant of the form $b_q$ is lower-bounded, for any $q$, by the inf-sup constant related to the hold-all domain $\hat{\Omega}$. Moreover, since the right-hand sides of \eqref{eq:stimauBrezzi} can be bounded by a data independent constant, also $\|\widetilde{\u}\|_{V_q}$, $\|\widetilde{p}\|_{P_q}$ are bounded, uniformly on $q$.\\
\label{oss:wellposedStilde}
\end{rmrk}

Finally, we introduce the cost functional
\begin{equation}
	J(q,\mathbf{u},p) = \int_{\Omega_q}{|\nabla \mathbf{u}|^2} + \alpha\|q''\|_{L^2(I)}^2 +\beta\left(\int_I {q(x)dx} - \overline{V}\right)^2,\\
\label{eq:minJ_NO}
\end{equation}
representing the total energy dissipation of the Stokes flow, with a regularization term $\|q''\|_{L^2(I)}^2$ (as in \cite{Kinigera}) and a volume penalty term, measuring the distance of the area under the graph of $q$ from a fixed value $\overline{V}$.
\footnote{Volume constraints are typical of shape optimization for fluid dynamics: see, e.g., \cite{MR0331973,Morin2012}.}\\

Let us introduce the state solution operator $\widetilde{S}(q)$, mapping each $q\in\Q$ to the corresponding solution $\widetilde{S}(q)=(\mathbf{u},p)$ of \eqref{eq:Stokesdeb}, and the reduced cost functional, as follows:
\beq
	\widetilde{j}:\Q\rightarrow\mathds{R},\quad \widetilde{j}(q)=J(q,\widetilde{S}(q)).
\label{eq:jtildedef}
\eeq

For convenience, it can be useful to define the following constants, whose existence is ensured by the fact that $q$ belongs to $\widehat{Q}^{ad}$:
\beq
	d_1,d_2>0 \text{ such that } \|q''\|_{L^\infty(I)}\leq d_1,\ |q'(0)|\leq d_2.
\label{eq:Bad}
\eeq

Finally, we introduce the set of admissible control variations, namely:
\beq
	\delta Q=\{\dq\in H^3(I)\cap H^1_0(I)\colon q+\dq\in\Q,\ \forall q\in \Q\}.
\label{eq:dQ}
\eeq

\begin{rmrk}
We point out that $\Q$ is convex, closed and bounded in $H^3(I)$: boundedness is stated in \eqref{eq:Qlim}, whereas closure and convexity are consequences of the fact that definitions \eqref{eq:Qad1} and \eqref{eq:Qlim} involve only constraints of the form $\zeta(q)\leq c$, where $c$ is a constant and $\zeta$ is a semi-norm in $H^3(I)$. Hence, closure follows from the continuity of any semi-norm in a Banach space, and convexity holds thanks to the triangle inequality.
\label{oss:Qconvesso}
\end{rmrk}

\subsection{Domain transformation}\label{sec:trasformazione}

In this section, we map the original problem \eqref{eq:Stokes} onto a reference-domain. The main advantage of this technique lays in the numerical solution of the optimization problem: solving the state problem on a reference domain avoids the need to deform the computational mesh at each step of the optimization algorithm.

Let us introduce the reference domain $\Omega_0=(0,1)^2$, which is equivalent to the choice $q\equiv 0$. It follows that any admissible domain $\Oq$ can be seen as a transformation of $\O$ by means of the map
\beq
	T_q :\Omega_0\rightarrow\Omega_q,\text{ \ with \ }T_q(x,y)=(\mathds{1}+V_q)(x,y)=\begin{pmatrix}x\\y+(1-y)q(x)\end{pmatrix}.
\label{eq:mapTq}
\eeq
We denote by $(\cdot,\cdot)$ the $L^2$ inner product on $\O$, while $(\cdot,\cdot)_I$ and $(\cdot,\cdot)_\Oq$ indicate the scalar product in $L^2(I)$ and $L^2(\Oq)$, respectively.\\

\begin{rmrk}[Notation I]
	We will use the following quantities depending on $T_q$:
	\begin{itemize}
		\item[] Map gradient: $DT_q$ with $(DT_q)_{i,j} = \partial_{x_j}(T_q)_i\qquad i,j=1,2$.
		\item[] Map jacobian: $\gamma_q = det(DT_q)$.
		\item[] Laplacian-related matrix: $A_q=\gamma_qDT_q^{-1}DT_q^{-T}$.
	\end{itemize}
\label{oss:coeffdef}
\end{rmrk}

\begin{rmrk}[Notation II]
By the superscript $\cdot^q$ we denote the composition with the map $T_q$. On the other hand, whenever no doubt arises on which $q$ is considered, the composition with the inverse map $T_q^{-1}$ will be denoted by $\widetilde{\cdot}$.
\end{rmrk}

We are now ready to state the variational problem \eqref{eq:Stokesdeb} on pulled-back spaces $V$ and $P$, that do not depend anymore on $q$:

\begin{pbdef}
	Find $(\mathbf{u},p)\in V\times P$, such that
	\beq
	\left\{
	\begin{aligned}
		a(q)(\mathbf{u},\mathbf{v}) + b(q)(\mathbf{v},p) &= F(q)(\mathbf{v}),\qquad&\forall\ \mathbf{v}\in  V,\\
		b(q)(\mathbf{u},\pi) &= G(q)(\pi),
			\qquad&\forall\ \pi \in P,\\
	\end{aligned}\right.
	\label{eq:StokesdebT}
	\eeq
	where
	\beq\begin{split}
		V&=\{\v\in[H^1(\O)]^2\colon \v=(v_x,v_y)=\mathbf 0\text{ on }\Gamma_3\cup\Gamma_0\text{ and }v_y=0\text{ on }\Gamma_2\},\\
		P&=L^2(\O),
	\end{split}\eeq
	and
	\beq\begin{split}
		a(q)(\mathbf{u},\mathbf{v})&=\int_{\Omega_0}{\left[\eta^q \mathbf{u}\cdot\mathbf{v}\gamma_q+\nu^q\ tr(\nabla \mathbf{u}A_q\nabla\mathbf{v}^T)\right]}d\Omega,\\
		b(q)(\mathbf{v},\pi) &= -\int_{\Omega_0}{\pi\ tr(\nabla\mathbf{v}DT_q^{-1})\gamma_q\,d\Omega},\\
		F(q)(\mathbf{v}) &= \int_{\Omega_0}{ \mathbf{f}^q\cdot \mathbf{v}\,\gamma_q\,d\Omega}-a(q)(\mathcal{R}\mathbf{g}_D,\mathbf{v}) + \int_{\Gamma_1}{\mathbf{g}_N\cdot\mathbf{v}\,d\Gamma},\\
		G(q)(\pi)&=-b(q)(\mathcal{R}\mathbf{g}_D,\mathbf{v}).
	\end{split}\eeq
\end{pbdef}
\begin{rmrk}[Lifting]
	$\mathcal{R}\mathbf g_D$ represents a continuous lifting of the Dirichlet datum $\mathbf g_D$ onto $\O$. However, as $\mathbf g_D$ is defined on $\Gamma_3$, where $T_q$ is equal to the identity, it does not need to be mapped onto the reference domain. In general, $\mathcal{R}\mathbf g_D\neq\widetilde{\mathcal{R}}\mathbf g_D\,\circ\,T_q$, but this is not a problem, since in the following we are not making use of any explicit expression of the lifting.
\end{rmrk}

Finally, we introduce the solution operator $S:\Q\to V\times P$, which maps an admissible control function to the solution of the transformed state problem \eqref{eq:StokesdebT}. It follows that the original optimization problem can be reformulated as follows:

\begin{pbdef}
	Find $\overline{q}\in\Q$ minimizing the functional $j$ defined in \eqref{eq:minj}, i.e.
	\beq
		j(\overline{q})=\min_{q\in\Q} j(q) =\min_{q\in\Q}J(q,S(q)\circ T_q^{-1}).
	\label{eq:minj}
	\eeq
\end{pbdef}
This is the formulation we will refer to on the rest of the paper.

\subsection{Well-posedness of the problem}\label{sec:wellposedness}

In this section, we analyze the well-posedness of the state problem \eqref{eq:StokesdebT} and the existence of an optimal solution to our minimization problem \eqref{eq:minj}.

At first, we observe that matrix $A_q$ 
belongs to $[L^\infty(\O)]^{2\times2}$, it is symmetric and positive definite, and its eigenvalues are lower-bounded by
\beq
	\overline{\lambda}=2\left(1+\frac{1+(d_1+d_2)^2}{\varepsilon} + \sqrt{\left(1+\frac{1+(d_1+d_2)^2}{\varepsilon}\right)^2-4}\right)^{-1}>0.
\label{eq:lambdamin}
\eeq
Under the same assumptions of Remark \ref{oss:wellposedStilde}, the coercivity of the form $a(q)$ and the continuity of the functionals and forms involved in \eqref{eq:StokesdebT} are given by the following inequalities, holding for any $\u,\v\in H^1_0(\O)$, $\pi\in L^2(\O)$, $q\in\Q$:
\begin{subequations}\begin{align}
	a(q)(\v,\v)&\geq\nu_0\overline{\lambda}\|\nabla\v\|^2=:\alpha_c\|\nabla\v\|^2, \label{eq:acoerc}\\
	\begin{split}
		|a(q)(\u,\v)|&\leq(\|\eta\|_{L^\infty(\hat{\Omega})}\|\gamma_q\|_\infty+\|\nu\|_{L^\infty(\hat{\Omega})}\|A_q\|_\infty)\|\nabla\u\|\|\nabla\v\|\leq\\
		&\leq\left(\|\eta\|_{L^\infty(\hat{\Omega})}(1+d_1+d_2)+\|\nu\|_{L^\infty(\hat{\Omega})}\frac{1}{\overline{\lambda}}\right)\|\nabla\u\|\|\nabla\v\|=:M\|\nabla\u\|\|\nabla\v\|,
	\end{split}\label{eq:acont}\\
	|b(q)(\v,\pi)|&\leq\|\gamma_qDT_q^{-T}\|_\infty\|\nabla\v\|\|\pi\|\leq (1+d_1+d_2)\|\nabla\v\|\|\pi\|=:M_b\|\nabla\v\|\|\pi\|, \label{eq:bcont}\\
	\begin{split}
	|F(q)(\v)|&\leq\|\gamma_q\|_\infty\|\mathbf f\|_{[L^2(\hat{\Omega})]^2}\|\v\|+Mc_\mathcal{R}\|\mathbf g_D\|_{H^{1/2}(\Gamma_3)}\|\nabla\v\|+\|\mathbf g_N\|_{[H^{-1/2}(\Gamma_1)]^2}c_{tr}\|\nabla\v\|\leq\\
	&\leq[c_{\hat{\Omega}}(1+d_1+d_2)\|\mathbf f\|_{[L^2(\hat{\Omega})]^2}+Mc_\mathcal{R}\|\mathbf g_D\|_{H^{1/2}(\Gamma_3)}+\|\mathbf g_N\|_{[H^{-1/2}(\Gamma_1)]^2}c_{tr}]\|\nabla\v\| =\\
	&=: M_F\|\nabla\v\|,
	\end{split}\label{eq:Fcont}
\end{align}\label{eq:contcoerc}\end{subequations}
where the constants $\alpha_c,M,M_b,M_F,M_G$ are independent of $q$.\\

The inf-sup condition for problem \eqref{eq:StokesdebT} reads


\begin{pbdef}
	There exists a positive constant $\overline{\beta}$, independent of $q$, such that
\beq
	\forall \pi \in P\ \ \exists\v\in V\ :\ b(q)(\v,\pi) \geq \overline{\beta} \|\nabla \v\|\|\pi\|.
\label{eq:infsup}
\eeq
\end{pbdef}
The validity of this property would allow to exploit the classical saddle-point-problem theory also for the transformed problem \eqref{eq:StokesdebT}.

To prove \eqref{eq:infsup}, we start considering the inf-sup condition on $\O$, with constant $\hat{\beta}>0$, namely:
\beq
	\forall \pi \in P\ \ \exists\v\in V \text{ such that } b(0)(\v,\pi)=-\int_\O\pi\,\div\,\v\,d\Omega\geq\hat{\beta}\|\pi\|\|\nabla\v\|.
\label{eq:infsup0}
\eeq
Employing the definition of $b(q)$, the following holds for any $q\in\Q$
\beq
	\begin{split}
	b(&q)(\v,\pi)=-\int_\O\pi\nabla\v\cdot\gamma_qDT_q^{-T}d\Omega = 
	-\int_\O\pi\nabla\v\cdot\left(\mathds{1}+\gamma_qDT_q^{-T}-\mathds{1}\right)d\Omega \geq\\
	&\geq -\int_\O \pi\,\div\,\v\,d\Omega-\left|\int_\O\pi\nabla\v\cdot\left(\gamma_qDT_q^{-T}-\mathds{1}\right)d\Omega\right|\geq(\hat{\beta}-\|\gamma_qDT_q^{-T}-\mathds{1}\|_\infty)\|\pi\|\|\nabla\v\|,
	\end{split}
\eeq
being $\pi\in P$ and $\v \in V$ related through \eqref{eq:infsup}.
As it holds
\beq
	\gamma_qDT_q^{-T}-\mathds{1}=\cof(DT_q)-\mathds{1}=\begin{pmatrix}-q & -(1-y)q'\\0&0\end{pmatrix},
\eeq
we get
\beq
	b(q)(\v,\pi)\geq(\hat{\beta}-\|q\|_{W^{1,\infty}(I)})\|\pi\|\|\nabla\v\|.
\eeq

Therefore, requiring $\|q\|_{W^{1,\infty}(\overline{I})}$ to be strictly smaller than $\hat{\beta}$ yields the validity of the inf-sup \eqref{eq:infsup}, uniformly on $q$.

It is easy to check (see, e.g., \cite{Duran2013}) that on the domain $\O$ the inf-sup constant $\hat{\beta}$ in \eqref{eq:infsup0} satisfies $\hat{\beta}\geq\frac{1}{4\sqrt{2}}$. Hence, in order to ensure the validity of \eqref{eq:infsup} it is sufficient to require
\beq
	\|q\|_{H^3(I)}\leq\frac{\xi}{4\sqrt{2}},\text{ \ for some }\xi\in(0,1),
	\label{eq:qxi}
\eeq
in the definition of the set $\Q$ of admissible controls.

\begin{rmrk}
	We remark that condition \eqref{eq:qxi} is representative of a class of sufficient conditions ensuring the validity of \eqref{eq:infsup}. Most likely, less stringent conditions can be found. However, real world shape optimization problems often deal with very smooth configurations, thus compatible with \eqref{eq:qxi}.
\end{rmrk}
Bearing in mind the properties showed at the beginning of this section, we can finally employ the classical results of saddle-point theory to prove the following result (see, e.g., \cite{Brezzi1974}):
\begin{prpstn}
	Under condition \eqref{eq:qxi}, for each $q\in\Q$ the pulled-back problem \eqref{eq:StokesdebT} admits unique solution, and the following inequality holds
	\beq
		\|S(q)\|_{V\times P}\leq c(\mathbf f,\mathbf g_D,\mathbf g_N,\eta,\nu,\hat{\Omega}),
	\eeq
	where the constant $c$ is independent of $q$.
\label{th:wellposedS}
\end{prpstn}

Concluding this section, we prove the existence of an optimal solution to \eqref{eq:minj}.
\begin{thrm}
	Let $\Q$ be a non-empty, convex, closed and bounded subset of $H^3(I)$ and let $S:\Q\to [H^1(\O)]^2\times L^2(\O)$ be the solution operator of problem \eqref{eq:StokesdebT}.
	Then, there exists a solution to the minimization problem \eqref{eq:minj}.
\label{th:esistenza}
\end{thrm}
\begin{proof}
	The proof follows standard ideas of calculus of variations. Hence, in the following we sketch the main steps of the proof.
	From Remark \ref{oss:Qconvesso}, we know that $\Q$ is a closed, bounded and convex subset of $H^3(I)$. This set is also non-empty, since $q\equiv 0$ fulfills all its constraints.\\
	Observing that $j(q)\geq 0$ for any $q\in\Q$ and that $\Q\neq\emptyset$, we have that a minimizing sequence $\{q_n\}_{n\in\mathds{N}}\subset\Q$ exists, such that
	\beq
		\lim_{n\in\mathds{N}}j(q_n) = \inf_{q\in Q^{ad}}j(q) =: \overline{j}.
	\eeq
	Being $\Q$ bounded in $H^3(I)$, the sequence $\{q_n\}$ is bounded itself, then there exists a subsequence $\{q_{n_k}\}$ and some $\overline{q}\in H^3(I)$ such that,
	\beq
		q_{n_k}\weak \overline{q}\quad \text{in }H^3(I)\qquad \text{for }k\to\infty.
	\eeq
	Being $\Q$ closed and convex, the limit $\overline{q}$ belongs to $\Q$.

	The next step to take is to show that we can take the limit also in the state variables sequence $\{S(q_{n_k})\}=\{(\u_k,p_k)\}$. For this purpose, following some ideas of the proof of Theorem 2.1 \cite{Gunzburger2000}, we consider the physical counterpart of the sequence, $\{\tilde{S}(q_{n_k})\}=\{S(q_{n_k})\circ T_q^{-1}\}$, and the trivial extension to zero of its elements in $\hat{\Omega}\supset\Oq$, denoted by $\{\hat{S}(q_{n_k})\}$.\\
	Thanks to the well-posedness of problem \eqref{eq:Stokes}, uniformly on $q\in\Q$, the sequence $\{\hat{S}(q_{n_k})\}$ is bounded in $\hat{V}\times\hat{P}=(H^1_0(\hat{\Omega})\cap V)\times (L^2(\hat{\Omega})\cap P)$. Hence, there exists a subsequence, for simplicity denoted by $\{S(q_l)\}$, and some $\overline{\hat{S}}=(\overline{\hat{\u}},\overline{\hat{p}})\in\hat{V}\times\hat{P}$ such that,
	\beq
		\widehat{S}(q_l)=(\widehat{\u}_l,\widehat{p}_l)\weak \overline{\widehat{S}}=(\overline{\widehat{\u}},\overline{\widehat{p}}) \text{ \ in \ }\widehat{V}\times \widehat{P}\qquad \text{for }l\to\infty.
	\eeq
	Now we have to prove that $\overline{S}=\overline{\hat{S}}|_{\Omega_{\overline{q}}}\circ T_{\overline{q}}$ is the state solution corresponding to $\overline{q}$, i.e. $\overline{S}=S(\overline{q})$. This can be done transforming each term in problem \eqref{eq:StokesdebT} back on $\Omega_{q_l}$, extending it on $\hat{\Omega}$ and then passing to the limit for $l\to\infty$. As a paradigmatic example, we consider the viscosity term. Taking $\widehat{\v}\in [C^\infty_0(\widehat{\Omega})]^2$, it holds
\beq
	\begin{split}
	\lim_{l\to\infty}&\int_\O\nu^{q_l}\nabla\u_lA_{q_l}\cdot\nabla\widehat{\v}\,d\Omega = \lim_{l\to\infty}\int_{\Omega_{q_l}}\nu\nabla\widetilde{\u}_l\cdot\nabla\widehat{\v}\,d\Omega =	\lim_{l\to\infty}\int_{\widehat{\Omega}}\nu\nabla\widehat{\u}_l\cdot\nabla\widehat{\v}\,d\Omega =\\
	 &= \int_{\widehat{\Omega}}\nu\nabla\overline{\widehat{\u}}\cdot\nabla\widehat{\v}\,d\Omega =
\int_\O\nu\nabla\overline{\u}A_{\overline{q}}\cdot\nabla\widehat{\v}\,d\Omega.
	\end{split}
\label{eq:limiti}
\eeq

	Finally, using dominate convergence theorem and the weak, lower semi-continuity of seminorms in a Banach space yields the weak, lower semi-continuity of functional $j$, allowing to conclude that
	\beq
		j(q_l)\to j(\overline{q})=\overline{j}\qquad \text{for }l\to\infty.
	\eeq
	Hence $\overline{q}$ turns out to be a solution of the optimization problem \eqref{eq:minj}.

\end{proof}


\subsection{Optimality conditions}\label{sec:ottimalita}

In this section, we inspect the first order optimality condition
\begin{equation}
	j'(\overline{q})(\dq)=0\qquad\forall\dq\in \delta Q,
\label{eq:ottimalita}
\end{equation}
in order to obtain the Hadamard formula (see, e.g., \cite{MR1179448}) for the gradient of functional $j$, useful for the analysis made in the following section and for numerical tests. 

We first recall the expression of $\widetilde{j}$ defined in \eqref{eq:jtildedef}, as
\beq
	\widetilde{j}(q) = \int_\Oq|\nabla\widetilde{\u}|^2\,d\Omega +\alpha\|q''\|_{L^2(I)}^2+\beta\left(\int_Iq(x)dx\ -\ \overline{V}\right)^2,
\label{eq:jtilde}
\eeq
where $\widetilde{\u}:\Oq\to\mathds{R}^2$, together with $\widetilde{p}:\Oq\to\mathds{R}$, is the solution of Stokes problem \eqref{eq:Stokes}.\\
The so-called shape-derivative of $(\widetilde{\u},\widetilde{p})$ can be defined as the solution $(\widetilde{\du},\widetilde{\delta p})$ of the following problem (see, e.g., \cite{Morin2012}):
\beq
	\left\{
	\begin{aligned}
	\eta\tilde{\du}-\div(\nu\nabla\tilde{\du})+\nabla\tilde{\delta p\,} &= 0, \qquad &\text{in } \Oq,\\
	\div(\tilde{\du}) &= 0\qquad &\text{in } \Oq,\\
	\nu\partial_{\mathbf n}\widetilde{\du}-\widetilde{\delta p}\mathbf n &= \mathbf 0, \qquad &\text{on } \Gamma_1,\\
	\partial_{\mathbf n}\widetilde{\delta u} =0,\quad\widetilde{\delta v}&=0,\qquad &\text{on } \Gamma_2,\\
	\widetilde{\du} &= \mathbf 0,\qquad &\text{on } \Gamma_3,\\
	\widetilde{\du} &=- (V_{q,\dq}\cdot\mathbf n)\partial_{\mathbf n}\widetilde{\u},\qquad &\text{on } \Gamma_q,\\
	\end{aligned}
	\right.
\label{eq:tildedu}
\eeq
where $V_{q,\dq}$ is the vector field describing a transformation from $\Oq$ to $\Omega_{q+\dq}$, given by
\beq
	V_{q,\dq}(x,y) = \begin{pmatrix}0\\ \frac{1-y}{1-q(x)}\dq(x)\end{pmatrix}.
\label{eq:Vqdq}
\eeq

Differentiating the expression \eqref{eq:jtilde} along direction $\dq$, one obtains
\beq
	\begin{split}
	j'(q)(\dq) &= 2(\nabla\widetilde{\u},\nabla\widetilde{\du})_\Oq + \int_{\Oq}{|\nabla\widetilde{\u}|^2V_{q,\dq}\cdot\mathbf n\,d\Gamma}+\\&+2\alpha(q'',\dq'')_I + 2\beta\left(\int_I{q(x)dx}-\overline{V}\right)\int_I{\dq(x)dx}.
	\end{split}
\label{eq:preHad}
\eeq

In order to make the dependence of $j'(q)(\dq)$ on $\dq$ completely explicit, we introduce the adjoint state $(\widetilde{z},\widetilde{s})$, solution of the following adjoint problem:
\begin{equation}
	\begin{cases}
		-\div(\nu\nabla\widetilde{\mathbf z}) +\eta\widetilde{\z}+\nabla \widetilde{s}= -2\Delta\widetilde{\u},\qquad &\text{in } \Oq,\\
		\div\,\widetilde{\mathbf z} = 0,\qquad &\text{in } \Oq,\\
		-\nu\partial_{\mathbf n}\widetilde{\z} +\widetilde{s}\mathbf n= -2\partial_{\mathbf n}\widetilde{\u},\qquad &\text{on } \Gamma_1,\\
		-\nu\partial_{\mathbf n}\widetilde{z}_x = -\nu\partial_{\mathbf n}\widetilde{z}_x + \widetilde{s}\,n_x = -2\partial_{\mathbf n} \widetilde{u}_x,\quad \widetilde{z}_y=0, &\text{on } \Gamma_2,\\
		\widetilde{\z} = \mathbf 0,\quad &\text{on } \Gamma_q\cup\Gamma_3.
	\end{cases}
\label{eq:PAHad}
\end{equation}

Using both problems \eqref{eq:tildedu} and \eqref{eq:PAHad}, and exploiting integration by parts and changes of variable from $\Oq$ to $\O$, and from $\Gamma_0$ to $I$, we can prove the following result (see, e.g., \cite{Fumagalli2013}):
\begin{lmm}
	Given the functional $j(q)$ defined as in \eqref{eq:minj}, its Gateaux-derivative in $q$ along direction $\dq$ is given by
\beq
	j'(q)(\dq) = 2\alpha(q'',\dq'')_I + (\Psi(q),\dq)_I\qquad \forall q\in\Q,\dq\in \delta Q,
\label{eq:HadPsi}
\eeq
where $\Psi(q):I\to\mathds{R}$ is defined as
\beq
	\begin{split}
	\Psi(q)(x)&=2\beta\left(\int_Iq(t)dt\ -\ \overline{V}\right)+ \\
	&\ +[\nabla\widetilde{\u}DT_q^{-1}DT_q^{-T}\mathbf n](x,q(x))\cdot[(\nu\nabla\z -\widetilde{\u})DT_q^{-1}DT_q^{-T}\mathbf n](x,q(x)).
	\end{split}
\label{eq:gradj}
\eeq
\label{th:Hadamard}
\end{lmm}

\section{A priori error estimates}\label{sec:stime}

In this section, we aim at deriving some a priori estimates for the numerical discretization error of the main quantities involved in our problem, namely the control function $q$, the state variable $S(q)$ and the reduced cost functional $j(q)$.

At first, we are going to discuss some differentiability properties of the state solution operator $S$, under suitable assumptions. Then, we will introduce a discretization on the control space and derive corresponding error estimates. Afterwards, the discretization of the state problem will be studied. Finally, we will derive a convergence result for the complete shape optimization problem.

\subsection{Solution operator properties}

In order to provide some differentiability properties for the state solution operator and the cost functional, we begin by considering the following generalization of the Implicit Function Theorem to Banach spaces:
\begin{thrm}[{\cite[Theorem~3.3]{Kinigera}}]
	Let $F\in C^k(X^{ad}\times Y,Z), k\geq 1$, where $Y$ and $Z$ are Banach spaces and $X^{ad}$ is an open subset of Banach space $X$. Suppose that $\mathcal{F}(x^*,y^*)=0$ and $\mathcal{F}'_y(x^*,y^*)$ is continuously invertible. Then there exist neighbourhoods $\Theta$ of $x^*$ in $X$, $\Phi$ of $y^*$ in $Y$ and a map $g\in C^k(\Theta,Y)$ such that $\mathcal{F}(x,g(y))=0$ for all $x\in\Theta$.
	Furthermore, $F(x,y)=0$ for $(x,y)\in\Theta\times\Phi$ implies $y=g(x)$.
\label{th:Dini}
\end{thrm}
As a direct consequence, we can prove the following result:
\begin{crllr}
	Let the following assumptions hold:
	\beq
		\eta,\nu\in C^2(\overline{\hat{\Omega}}),\qquad\mathbf f\in[C^2(\overline{\hat{\Omega}})]^2.
	\eeq
	Then, the solution operator $S$ is at least twice continuously Fréchet-differentiable.
\label{th:Dinicor}
\end{crllr}
\begin{proof}
	It is enough to use Theorem \ref{th:Dini}, with $X=H^2(I)\cap H^1_0(I), Y=V\times P, Z=Y^*$, the open set $X^{ad}=int(\Q)$ and the map $\mathcal{F}:X^{ad}\times Y \to Z$ such that
	\beq
	\mathcal{F}(q;\u,p)=
	\begin{pmatrix}
		a(q)(\u,\cdot)+b(q)(\cdot,p)-F(q)(\cdot)\\
		b(q)(\u,\cdot)-G(q)(\cdot)
	\end{pmatrix}\text{ for any }q\in int(\Q), \u\in V, p \in P.
	\eeq
	The regularity of the map is a consequence of the regularity of the forms involved in its definition.\\
	It is easy to check that the operator $S$ corresponds to the map $g$ defined in Theorem \ref{th:Dini}, hence the regularity result for $g$ holds for $S$ as well.
\end{proof}

Now, let us preliminarily collect some properties of the map $T_q$.
\begin{prpstn}
	Given $q\in\Q$, the maps defined in Remark \ref{oss:coeffdef}, depending on $T_q$ and its derivatives, satisfy the following inequalities, for any admissible variation $\dq\in \delta Q$:
	\begin{enumerate}
	\item $\|\gamma'_{q,\dq}\|_\infty=\|\div(V_\dq)\|_\infty\leq c\|\dq\|_{L^\infty(I)}\leq \bar{c}\|\dq\|_{H^1(I)},$
	\item $\|V_\dq\|_\infty\leq c\|\dq\|_{H^1(I)},$
	\item $\|\cof(DV_\dq)\|_\infty=\|DV_\dq\|_\infty\leq c\|\dq\|_{H^2(I)},$
	\item $\div(\cof(DV_\dq)) = (0,0)^T,$
	\item $\|A'_{q,\dq}\|_\infty\leq c\|\dq\|_{H^2(I)},$
	\item $\|\div(A'_{q,\dq})\|_\infty\leq c\|\dq\|_{H^2(I)},$
	\item $\|A'_{q,\dq}\|_2\leq c\|\dq\|_{H^1(I)},$
	\item $\|A''_{q,\dq,\dq}\|_\infty\leq c\|\dq\|_{H^2(I)}^2,$
	\end{enumerate}
	where the constants $c$ and $\overline{c}$ are independent of $q$ and $\dq$.
\label{th:coeff}
\end{prpstn}
\begin{proof}
	The result simply follows from direct computation and the application of the Fundamental Theorem of Calculus.
\end{proof}

The differentiability properties of the solution operator $S$ are characterized in the following result.
\begin{thrm}
\label{th:dS}
	The first and second variations of the solution operator $S$ along the directions $\delta q, \tau q \in Q^{ad}$ are defined as follows:
	\begin{enumerate}
		\item $S'(q)(\delta q)=(\delta \u, \delta p) \in V\times P$, where $(\delta \u, \delta p)$ is the solution of
			\begin{equation}
			\left\{
			\begin{aligned}
			a(q)(\du,\mathbf{v}) + b(q)(\mathbf{v},\delta p) =\ &\dot{F}(q,\dq)(\mathbf{v}) - \dot{a}(q,\dq)(\u,\mathbf v) - \dot{b}(q,\dq)(\mathbf v, p)\qquad&\forall\ \mathbf{v}\in  V,\\
			b(q)(\du,\pi) =\ &\dot{G}(q,\dq)(\pi)-\dot{b}(q,\dq)(\u,\pi)\qquad&\forall\ \pi \in P.\\
			\end{aligned}\right.\\
			\label{eq:dS}
			\end{equation}
		\item $S''(q)(\dq,\tau q)=(\tau\du,\tau\delta p) \in V\times P$, where $(\tau\du,\tau\delta p)$ is the solution of
			\begin{equation}
			\left\{
			\begin{aligned}
			&\begin{aligned}
			a(q)&(\tau\du,\mathbf{v}) + b(q)(\mathbf{v},\tau \delta p) = \\
			&=\ddot{F}(q,\dq,\tau q)(\mathbf{v}) - \ddot{a}(q,\dq,\tau q)(\u,\mathbf v) - \ddot{b}(q,\dq,\tau q)(\mathbf v, p)+\\
			&- \dot{a}(q,\dq)(\tu,\mathbf v) - \dot{b}(q,\dq)(\mathbf v, \tau p) - \dot{a}(q,\tau q)(\du,\mathbf v) - \dot{b}(q,\tau q)(\mathbf v, \delta p)
			\end{aligned}\quad&\forall\ \mathbf{v}\in  V,\\
			&\begin{aligned}
			b(q)&(\tau\du,\pi) = \ddot{G}(q,\dq,\tau q)(\pi)-\ddot{b}(q,\dq,\tau q)(\u,\pi) +\\
			&- \dot{b}(q,\dq)(\tu,\pi) - \dot{b}(q,\tau q)(\du,\pi)
			\end{aligned}\quad&\forall\ \pi \in P,\\
			\end{aligned}\right.
			\label{eq:d2S}
			\end{equation}
			with $(\tu,\tau p) = S'(q)(\tau q)$.
	\end{enumerate}
	The forms and functionals employed in \eqref{eq:dS} and \eqref{eq:d2S} are defined as follows:
	\begin{equation}
	\begin{split}
	\dot{F}(q,\dq)(\mathbf v) &= \int_{\Omega_0}\left(\gamma'_{q,\dq}\mathbf{f}^q\cdot\mathbf v + \gamma_q\nabla\mathbf f^q V_\dq\cdot\mathbf v\right)d\Omega-\dot{a}(q,\dq)(\mathcal{R}\mathbf g_D,\mathbf v),\\
	\dot{G}(q,\dq)(\pi) &= -\dot{b}(q,\dq)(\mathcal{R}\mathbf g_D,\pi),\\
	\dot{a}(q,\dq)(\u,\mathbf v) &= \int_\O\left[\left(\gamma_q\nabla\eta^q\cdot V_\dq + \eta^q\gamma'_{q,\dq}\right)\u\cdot\mathbf v \right. +\\
			&+\left. \nabla\nu^q\cdot V_\dq tr(\nabla\u A_q\nabla\mathbf v^T) + \nu^q tr(\nabla\u A'_{q,\dq}\nabla\v^T)\right]d\Omega,\\
	\dot{b}(q,\dq)(\mathbf v,\pi) &= -\int_\O\pi\nabla\v\cdot \cof(DV_\dq)\,d\Omega,\\
	\ddot{F}(q,\dq,\tau q)(\mathbf{v}) &= \int_\O[\gamma''_{q,\dq\tau q}\mathbf{f}^q\cdot\mathbf v + \gamma'_{q,\dq}\nabla\mathbf f^q V_{\tau q}\cdot\mathbf v + \gamma'_{q,\tau q}\nabla\mathbf f^q V_\dq\cdot\mathbf v +\\
			&+ \gamma_q(\widetilde{\nabla}^2\mathbf f^q V_{\tau q} + \nabla\mathbf f^qDV_{\tau q})V_\dq\,\cdot\v]d\Omega-\ddot{a}(q,\dq)(\mathcal{R}\mathbf g_D,\mathbf v),\\
	\ddot{G}(q,\dq,\tau q)(\pi) &= -\ddot{b}(q,\dq,\tau q)(\mathbf v, \pi) = 0,\\
	\ddot{a}(q,\dq,\tau q)(\u,\mathbf v) &= \int_\O\{[\gamma'_{q,\tau q}\nabla\eta^q\cdot V_\dq + (\nabla^2\eta^q V_{\tau q} + DV_{\tau q}^T\nabla\eta^q)\cdot V_\dq\gamma_q + \\
			&+ \eta^q\gamma''_{q,\dq,\tau q})\u\cdot\mathbf v + \gamma'_{q,\dq}\nabla\eta^q\cdot V_{\tau q}]\u\cdot\v +\\
			&+ (\nabla^2\nu^q\,V_\dq + DV_{\tau q}^T\nabla\nu^q)\cdot V_\dq\, tr(\nabla\u A_q\nabla\mathbf v^T) +\\
			&+ \nabla\nu^q\cdot V_\dq\, tr(\nabla\u A'_{q,\tau q}\nabla\mathbf v^T) + \nu^q\, tr(\nabla\u A''_{q,\dq,\tau q}\nabla\v^T)+\\
			&+ \nabla\nu^q\cdot V_{\tau q}\,tr(\nabla\u A'_{q,\dq}\nabla\v^T)\}d\Omega,\\
	\ddot{b}(q,\dq,\tau q)(\mathbf v, \pi) &= 0,
	\end{split}
	\label{eq:formedot}	
	\end{equation}
	with the differential operator $\widetilde{\nabla}^2$ acting as $\left(\widetilde{\nabla}^2\boldsymbol\varphi\right)_{ijk} = \left(\nabla^2\varphi_i\right)_{kj}$ and  the over-signed dots denoting the partial Gateaux derivative w.r.t.\ the control $q$.
	Moreover, the following stability results hold:
	\beq
		\begin{split}
		&\|S(q)\|_{V\times P}\leq c,\\
		&\|S'(q)(\dq)\|_{V\times P}\leq c\|\dq\|_{H^2(I)},\\
		&\|S''(q)(\dq,\dq)\|_{V\times P}\leq c\|\dq\|_{H^2(I)}^2,
		\end{split}
	\label{eq:dS_stab}
	\eeq
	provided that the data satisfy the following regularity requirements:
	\beq
		\eta,\nu\in W^{2,\infty}(\hat{\Omega}),\quad \mathbf f\in [H^1(\hat{\Omega})]^2.
	\label{eq:dS_datareg}
	\eeq
\end{thrm}
\begin{proof}
	The weak problems defined in \eqref{eq:dS}-\eqref{eq:d2S} can directly be obtained by differentiating the state problem \eqref{eq:StokesdebT} w.r.t.\ $q$.
The stability results \eqref{eq:dS_stab} follow from classical well-posedness results for saddle-point problems (see, e.g., \cite{Girault1986}), combined with  Proposition \ref{th:coeff} (see \cite{Fumagalli2013} for details).
\end{proof}
\begin{rmrk}
We observe that the first derivative of the solution operator, $S'(q)(\dq)=(\du,\delta p)$, is the transformation of the shape derivative $(\tilde{\du},\tilde{\delta p\,})$ introduced in \eqref{eq:tildedu}, since one can prove that $\du=\tilde{\du}\circ T_q,\ \delta p=\tilde{\delta p\,}\circ T_q$.
\end{rmrk}

Hinging upon Theorem \ref{th:dS}, we are now ready to compute the derivatives of $j$, as follows:
\begin{subequations}
\label{eq:jj'j''}
\begin{align}
	j(q)&=(\nabla\u\, A_q,\nabla\u) + \alpha\|q''\|_I^2 + \beta\left(\int_Iq(x)dx\ - \ \overline{V}\right)^2 \label{eq:j},\\
	\begin{split} 
		j'(q)(\dq) &= (\nabla\u\, A'_{q,\dq},\nabla\u) + 2(\nabla\du\, A_q,\nabla\u) + 2\alpha(\dq'',q'')_I +\\
		&+ 2\beta\left(\int_I{q(x)dx}-\overline{V}\right)\int_I{\dq(x)dx},
	\end{split}\label{eq:j'}\\
	\begin{split}
		j''(q)(\dq,\tau q) &= (\nabla\u\, A''_{q,\dq,\tau_q},\nabla\u) + 2(\nabla\tau\u\, A'_{q,\dq}+\nabla\du\, A'_{q,\tau q},\nabla\u) +\\
		&+ 2(\nabla\du\,A_q,\nabla\tau\u) + 2(\nabla\tau\du\,A_q,\nabla\u)+\\
		&+ 2\alpha(\dq'',\tau q'')_I + 2\beta\int_I{\dq(x)dx}\int_I{\tau q(x)dx},
	\end{split}\label{eq:j''}
\end{align}
\end{subequations}
where $\u,\du,\tau\u, \tau\du$ are the same as in Theorem \ref{th:dS}. The continuity of the derivatives is an easy consequence of the regularity and symmetry of the matrix $A_q$ and its derivatives.

\subsection{Control discretization}

Let $\{I_i=(x_{i-1},x_i)\}_{i=1}^N$ be a partition of the domain $I$, with discretization parameter $\sigma=\max_{i\in\{1,\dots,N\}}|I_i|$. We can then define the discrete controls set as
\beq
	\Q_\sigma = \Q\cap Q_\sigma,\text{  with  }Q_\sigma=\{q\in C^0(\overline{I})\colon q|_{I_i}\in\mathds{P}_4(I_i),\ i\in\{1,\dots,N\}\}.
\eeq
The semi-discretized optimization problem reads as follows
\beq
	\min_{q_\sigma\in\Q_\sigma} j(q_\sigma) = J(q_\sigma,\widetilde{S}(q_\sigma)).
\label{eq:minJsigma}
\eeq
As $\Q\supseteq\Q_\sigma$, the minimization problem \eqref{eq:minJsigma} inherits the existence and regularity properties holding for the original continuous optimization problem \eqref{eq:minj}.

Let us denote by $\Pi^4_\sigma\colon L^2(I)\to Q_\sigma$ the classical polynomial interpolation operator and notice that $\Pi^4_\sigma(\Q)\subseteq\Q_\sigma$. Standard interpolation error estimates hold (see e.g. \cite{Brenner2008}): for $r\geq1$, $0\leq m\leq r+1$, it holds that
\beq
	|q-\Pi^r_\sigma q|_{H^m(I)}\leq c\sigma^{r+1-m}|q|_{H^{r+1}(I)}\quad\forall q\in H^{r+1}(I).
\label{eq:qinterpr}
\eeq
In this section, we aim at proving the following convergence result:
\begin{prpstn}
	Let $\overline{q}\in\Q$ be the exact solution of \eqref{eq:minj}, and $\overline{q}_\sigma$ the solution of the partially discretized problem \eqref{eq:minJsigma}. Then, assuming that the optimal control $\overline{q}$ belongs to $H^5(I)$, the following convergence error estimate holds:
	\beq
		\|\overline{q}-\overline{q}_\sigma\|_{H^3(I)}\leq c\sigma^2|\overline{q}|_{H^5(I)}.
	\label{eq:stimaerrq}
	\eeq
\label{th:stimaerrq}
\end{prpstn}
\begin{rmrk}
	We observe that Proposition \ref{th:stimaerrq} needs the optimal control $\overline{q}$ to be in $H^5(I)$. To achieve this regularity, there is no need to re-define the admissible controls set $\Q$, but it is sufficient to assume the validity of a regularity result for the classical Stokes problem. This assumption and the proof of the needed regularity on $\overline{q}$ are reported in Appendix \ref{app:regolarita}.
\end{rmrk}
In order to prove Proposition \ref{th:stimaerrq}, we need to collect some preliminary results that will be derived under the following two assumptions, already employed in \cite{Kinigera}.

\begin{ssmptn}[{\cite[Assumption 1.5]{Kinigera}}]
	For the optimal solution $\overline{q}$ of problem \eqref{eq:minj}, the constraint $q\leq1-\varepsilon$ is not active, i.e.
	\beq
		\exists\delta>0\text{ such that }\overline{q}(x)\leq1-\varepsilon-\delta\quad\forall x\in I.
	\eeq
\label{ass:noeps}
\end{ssmptn}
\begin{ssmptn}[{\cite[Assumption 3.1]{Kinigera}}]
	For any local minimum $\overline{q}$, we have
	\beq
		j''(\overline{q})(\dq,\dq)>0\quad\forall\dq\in \delta Q\setminus\{0\}.
	\eeq
\label{ass:jpos}
\end{ssmptn}

We start by proving some regularity results for the solution operator $S$ and its derivatives.
\begin{lmm}
	Let $S$ be the solution operator of the transformed Stokes problem \eqref{eq:StokesdebT}. If there exists some $k>0$ such that data functions fulfill the regularity requests
	\beq
		\eta,\nu\in C^k(\hat{\Omega}),\qquad \mathbf f\in [C^k(\hat{\Omega})]^2,
	\eeq
	then $S$ is at least $k$ times continuously Fréchet differentiable.
\label{th:SCk}
\end{lmm}
\begin{proof}
	The proof is the same as in Corollary \ref{th:Dinicor}, simply applying the Implicit Function Theorem in the form presented in Theorem \ref{th:Dini}.
\end{proof}
Based on the previous result, we can prove the following:
\begin{lmm}
	Let $k\in\mathds{N}$ and let data functions fulfill the following regularity requests:
	\beq
		\eta,\nu\in C^{k+1}(\hat{\Omega}),\qquad \mathbf f\in [C^{k+1}(\hat{\Omega})]^2.
	\eeq
	Then, for any $q,r\in \Q$ and $\dq_1,\dq_2,\dots,\dq_k\in \delta Q$, the following inequalities hold:
	\beq
		\|S^{(i)}(q)(\dq_1,\dots,\dq_i)-S^{(i)}(r)(\dq_1,\dots,\dq_i)\|_{V\times P}\leq c\|q-r\|_{H^2(I)}\prod_{j=1}^{i}\|\dq_j\|_{H^2(I)},\quad \text{for }i=0,\dots,k.
	\eeq
\label{th:SLip}
\end{lmm}
\begin{proof}
Let $q$ and $r$ be two control functions in $\Q$ and $\delta q,\tau q$ admissible control variations. Applying Lemma \ref{th:SCk} under the hypotheses of the present Lemma, we get
\beq
	S\in C^{k+1}(int(Q^{ad}); V\times P).
\eeq
Let us consider $k=0$. As $S\in C^1$, given the control functions $q,r\in\Q$, the Mean Value Theorem ensures that
\beq
	\exists \xi\in\Q\text{ such that } S(q)-S(r) = S'(\xi)(q-r).
\eeq
Being the Fréchet derivative $S'(\xi)$ a linear operator on the control variation, its continuity is equivalent to its boundedness, thus we get
\beq
	\|S(q)-S(r)\|_{V\times P}=\|S'(\xi)(q-r)\|_{V\times P}\leq c\|q-r\|_{H^2(I)}.
\eeq

In the general case $k>0$, for each $i\in\{0,\dots,k\}$ there exists $\xi_i\in \Q$ such that
\beq
	S^{(i)}(q)-S^{(i)}(r)= S^{(i+1)}(\xi_i)(q-r),
\label{eq:TaylordS}
\eeq
where we remark that \eqref{eq:TaylordS} is an equality between linear operators belonging to $\mathscr{L}_i:=\mathscr{L}(\delta Q^i;V\times P)$. Observing that $S^{(i+1)}(\xi_i)\in\mathscr{L}_{i+1}$, we can proceed as before to obtain
\beq
	\begin{split}
	\|S^{(i)}(q)(&\delta q_1,\dots,\dq_{i})-S^{(i)}(r)(\dq_1,\dots,\dq_i)\|_{V\times P}=\\
	& = \|S^{(i)}(q)-S^{(i)}(r)\|_{\mathscr{L}_i} \prod_{j=1}^{i}\|\delta q_i\|_{H^2(I)}= \|S^{(i+1)}(\xi_i)(q-r)\|_{\mathscr{L}_{i}} \prod_{j=1}^{i}\|\delta q_i\|_{H^2(I)}\leq\\
	& \leq  \|S^{(i+1)}(\xi_i)\|_{\mathscr{L}_{i+1}} \|q-r\|_{H^2(I)}\prod_{j=1}^{i}\|\delta q_i\|_{H^2(I)}.
	\end{split}
\eeq
Since $S^{(i+1)}$ is continuous, it is also bounded, so there exists a constant $c>0$ such that \linebreak[4]$\|S^{(i+1)}(\xi)\|_{\mathcal{L}_{i+1}}\leq c$ for all $\xi\in\delta Q$. Hence the proof is complete.

\end{proof}

The continuity of the solution operator $S$ directly implies the continuity of the functional $j$, as stated in the following result.
\begin{lmm}
	For any $q,r\in \Q$ and any $\dq\in H^2(I)\cap H^1_0(I)$, it holds that
	\begin{enumerate}[(a)]
	\item $|j(q)-j(r)|\leq c\|q-r\|_{H^2(I)},$
	\item $|j'(q)(\dq)-j'(r)(\dq)|\leq c\|q-r\|_{H^2(I)}\|\dq\|_{H^2(I)},$
	\item $|j''(q)(\dq,\dq)-j''(r)(\dq,\dq)|\leq c\|q-r\|_{H^2(I)}\|\dq\|_{H^2(I)}^2.$
	\end{enumerate}
\label{th:jLip}
\end{lmm}
\begin{proof}
	Let us fix $q,r\in\Q$. To simplify the notation, let $S(q)=(\u,p)$ and $S(r)=(\z,s)$.
	As the proofs of (a)-(c) are similar, we focus on (c), highlighting the most technical parts.\\
	Bearing in mind the expression of $j''$ (see \eqref{eq:j''}), we first focus on the following term:\noeqref{eq:j}\noeqref{eq:j'}
	\begin{equation*}\begin{split}
	&|\left(\nabla\u A''_{q,\dq,\dq},\nabla\u\right)-\left(\nabla\z A''_{r,\dq,\dq},\nabla\z\right)|=\\
	&=|\left((\nabla\u-\nabla\z)A''_{q,\dq,\dq},\nabla\u+\nabla\z\right)+\left(\nabla\z(A''_{q,\dq,\dq}-A''_{r,\dq,\dq}),\nabla\z\right)|\leq\\
	&\leq\|A''_{q,\dq,\dq}\|_\infty\left(\|\nabla\u\|+\|\nabla\z\|\right)\|\nabla\u-\nabla\z\|+\|\nabla\z\|^2\|A''_{q,\dq,\dq}-A''_{r,\dq,\dq}\|_\infty\leq\\
	&\leq c\|\dq\|_{H^2(I)}^2,
	\end{split}\end{equation*}
	where the state variables have been bounded using Lemmas \ref{th:wellposedS} and \ref{th:dS}, while $\|A''_{q,\dq,\dq}\|_\infty$ has been handed employing Proposition \ref{th:coeff}.\\
	Using the same results, it is easy to bound also the following term:
	\begin{equation*}\begin{split}
	&|\left(\nabla\du\,A'_{q,\dq},\nabla\u\right)-\left(\nabla\delta\z\,A'_{r_\dq},\nabla\z\right)|=\\
	&=|\left((\nabla\du-\nabla\delta\z)A'_{q,\dq},\nabla\u\right)+\left(\nabla\delta\z(A'_{q,\dq}-A'_{r,\dq}),\nabla\u\right)+\left(\nabla\delta\z\, A'_{r,\dq},\nabla\u-\nabla\z\right)|.
	\end{split}\end{equation*}
	All the other terms entering in $j''$ can be treated in a similar way, to get (c).
\end{proof}
%

The results stated so far are sufficient to prove the following coercivity result on $j$.
\begin{lmm}
	If $\overline{q}$ is a local solution of \eqref{eq:minj}, fulfilling Assumption \ref{ass:jpos}, then there exist $\delta_1,\delta_2>0$ such that, if $\|\overline{q}-r\|_{H^2(I)}\leq\delta_1$ for $r\in\Q$, then
	\beq
		j''(r)(\dq,\dq)\geq\frac{\delta_2}{2}\|\dq\|^2_{H^3(I)}\qquad\forall\dq\in\Q.
	\eeq
\label{th:314}
\end{lmm}
\noindent The proof of Lemma \ref{th:314} is the same as in \cite[Lemma 3.14]{Kinigera}, replacing $H^2(I)$ with $H^3(I)$ and provided two more intermediate results, reported in Appendix \ref{sec:per314}.\\


Now, we are ready to conclude this section with the proof of Proposition \ref{th:stimaerrq}.
\begin{proof}[Proof (Proposition \ref{th:stimaerrq})]
The Mean Value Theorem and Lemma \ref{th:314} imply the existence of some $t\in[0,1]$ such that, for $\xi=t\,\Pi^4_\sigma\overline{q}+(1-t)\overline{q}_\sigma$, we have
\beq
	\begin{split}
	\frac{\delta_2}{2}\|\Pi^4_\sigma\overline{q}-\overline{q}_\sigma\|_{H^3(I)}^2 &\leq j''(\xi)(\Pi^4_\sigma\overline{q}-\overline{q}_\sigma,\Pi^4_\sigma\overline{q}-\overline{q}_\sigma) =\\
	&= j'(\Pi^4_\sigma\overline{q})(\Pi^4_\sigma\overline{q}-\overline{q}_\sigma)-j'(\overline{q}_\sigma)(\Pi^4_\sigma\overline{q}-\overline{q}_\sigma) =\\
	&\overset{\text{a}}{=} j'(\Pi^4_\sigma\overline{q})(\Pi^4_\sigma\overline{q}-\overline{q}_\sigma)-j'(\overline{q})(\Pi^4_\sigma\overline{q}-\overline{q}_\sigma) \leq\\
	&\overset{\text{b}}{\leq}  c\|\overline{q}-\Pi^4_\sigma\overline{q}\|_{H^3(I)}\|\Pi^4_\sigma\overline{q}-\overline{q}_\sigma\|_{H^3(I)}\leq\\
	&\overset{\text{c}}{\leq} c\,\sigma^2\|\overline{q}\|_{H^5(I)}\|\Pi^4_\sigma\overline{q}-\overline{q}_\sigma\|_{H^3(I)},
	\end{split}
\label{eq:fineerrq}
\eeq
where we used:
\begin{enumerate}[\quad(a)]
	\item $j'(\overline{q})(\Pi^4_\sigma\overline{q}-\overline{q}_\sigma)=j'(\overline{q}_\sigma)(\Pi^4_\sigma\overline{q}-\overline{q}_\sigma) = 0$, due to Assumption \ref{ass:noeps} and then the first order optimal condition;
	\item point \emph{b} of Lemma \ref{th:jLip} and the fact that $\|\cdot\|_{H^2(I)}\leq\|\cdot\|_{H^3(I)}$;
	\item the interpolation error estimate \eqref{eq:qinterpr}.
\end{enumerate}
From \eqref{eq:fineerrq}, we obtain
\beq
	\|\Pi_\sigma\overline{q}-\overline{q}_\sigma\|_{H^3(I)}\leq c\sigma^2\|\overline{q}\|_{H^5(I)}.
\eeq
Finally, triangular inequality gives the thesis.
\end{proof}

\subsection{State discretization}\label{sec:state_discr}

Let $\mathcal{T}_h$ be a regular triangulation of $\O$, with discretization parameter $h=\smash{\displaystyle\max_{K\in\mathcal{T}_h}}|K|$. We can thus introduce the finite element spaces
\begin{equation}\begin{split}
	X_h^r(\O) &=\{\varphi\in C^0(\overline{\Omega}_0)\ \colon\ \varphi|_K\in\mathds{P}_r(K)\ \ \forall K\in\mathcal{T}_h\},\\
	V_h &= V\cap[X_h^2(\O)]^2,\\
	P_h &= P\cap X_h^1(\O),
\label{eq:spacesh}
\end{split}\end{equation}
where $\mathds{P}_r(K)$ is the space of polynomials on $K$ having degree less than or equal to $r$.

Passing from the continuous to the discrete case, the variational forms involved in problem \eqref{eq:StokesdebT} preserve all their properties, with discrete inf-sup condition ensured by the following:
\begin{prpstn}[LBB condition]
There exists a positive constant $\overline{\beta}$ such that
\beq
	\forall \pi_h \in P_h\ \exists\v_h\in V_h\ :\ b(q)(\v_h,\pi_h) \geq \overline{\beta} \|\nabla \v_h\|\|\pi_h\|,
\label{eq:LBB}
\eeq
and $\overline{\beta}$ is independent from $q\in\Q$ and from $h\in[0,\widehat{h}]$, for a certain $\hat{h}>0$.
\label{th:LBB}
\end{prpstn}
\begin{proof}
	From FEM approximation of Stokes problem \cite{Girault1986}, we know that pair $(V_h,P_h)$ is stable, i.e. there exists a constant $\widehat{\beta}>0$ such that
	\beq
		\forall \pi_h \in P_h\ \exists\v_h\in V_h\ :\ b(0)(\v_h,\pi_h) \geq \widehat{\beta} \|\nabla \v_h\|\|\pi_h\|,
	\label{eq:LBB0}
	\eeq
	with $\hat{\beta}$ independent from $h\in[0,\hat{h}]$.

	 In order to show that such discrete spaces fulfill inf-sup condition also for the transformed form $b(q)$, one can just follow the steps presented in section \ref{sec:wellposedness}, with constant $\hat{\beta}$ from \eqref{eq:LBB0}. Indeed, no assumptions on the spaces $V,P$ have been made there, apart from the validity of inf-sup condition for $b(0)$.
	\label{oss:LBB}
\end{proof}

The finite element discretization of \eqref{eq:StokesdebT} reads as follows:
\begin{equation}
\begin{split}
	&\text{Find }(\mathbf{u}_h,p_h)\in V_h\times P_h, \text{ such that } \\
	&\left\{
	\begin{aligned}
		a(q_\sigma)(\mathbf{u}_h,\mathbf{v}_h) + b(q_\sigma)(\mathbf{v}_h,p_h) &= F(q_\sigma)(\mathbf{v}_h)\qquad&\forall\ \mathbf{v}_h\in  V_h,\\
		b(q_\sigma)(\mathbf{u}_h,\pi_h) &= G(q_\sigma)(\pi_h) \qquad&\forall\ \pi_h \in P_h.\\
	\end{aligned}\right.\\
\end{split}
\label{eq:StokesdebTh}
\end{equation}
The well-posedness of \eqref{eq:StokesdebTh} stems from the validity of \eqref{eq:LBB}.

The discrete state solution operator, resulting from problem \eqref{eq:StokesdebTh}, and the corresponding discrete cost functional, are defined as
\beq
	S_h:\Q\to V_h\times P_h, \text{ with } S_h(q)=(\u_h,p_h),\qquad j_h\colon\Q\to\mathds{R}, \text{ with } j(q)=J(q,S_h(q)\circ T_q^{-1}),
\eeq
whereas the fully discretized shape optimization problem can be written as
\beq
	\min_{q_\sigma\in\Q_\sigma} j_h(q_\sigma)=J(q_\sigma,S_h(q_\sigma)\circ T_{q_\sigma}^{-1}).
\label{eq:minJh}
\eeq

\noindent For future use, it is useful to explicitly write the problems defining the derivatives of $S_h$:
	\begin{enumerate}
		\item $S_h'(q)(\delta q)=(\delta \u_h, \delta p_h) \in V_h\times P_h$, where $(\delta \u_h, \delta p_h)$ is the solution of
			\begin{equation}
			\left\{
			\begin{aligned}
			a(q)(\du_h,\mathbf{v}_h) + &b(q)(\mathbf{v}_h,\delta p_h) =\dot{F}(q,\dq)(\mathbf{v_h}) + &\\
			&\quad- \dot{a}(q,\dq)(\u_h,\mathbf v_h) - \dot{b}(q,\dq)(\mathbf v_h, p_h)\qquad&\forall\ \mathbf{v_h}\in  V_h,\\
			b(q)(\du_h,\pi_h) =\ &\dot{G}(q,\dq)(\pi_h)-\dot{b}(q,\dq)(\u_h,\pi_h)\qquad&\forall\ \pi_h \in P_h.\\
			\end{aligned}\right.\\
			\label{eq:dSh}
			\end{equation}
		\item $S_h''(q)(\dq,\tau q)=(\tau\du_h,\tau\delta p_h) \in V_h\times P_h$, where $(\tau\du_h,\tau\delta p_h)$ is the solution of
			\begin{equation}
			\left\{
			\begin{aligned}
			a(q)&(\tau\du_h,\mathbf{v}_h) + b(q)(\mathbf{v}_h,\tau \delta p_h) = \\
			&=\ddot{F}(q,\dq,\tau q)(\mathbf{v}_h) - \ddot{a}(q,\dq,\tau q)(\u_h,\mathbf v_h) - \ddot{b}(q,\dq,\tau q)(\mathbf v_h, p_h)+\\
			&\qquad- \dot{a}(q,\dq)(\tu_h,\mathbf v_h) - \dot{b}(q,\dq)(\mathbf v_h, \tau p_h) + \\
			&\qquad- \dot{a}(q,\tau q)(\du_h,\mathbf v_h) - \dot{b}(q,\tau q)(\mathbf v_h, \delta p_h)
			\quad&\forall\ \mathbf{v_h}\in  V_h,\\
			b(q)&(\tau\du_h,\pi_h) = \ddot{G}(q,\dq,\tau q)(\pi_h)-\ddot{b}(q,\dq,\tau q)(\u_h,\pi_h) +\\
			&\qquad\qquad\qquad- \dot{b}(q,\dq)(\tu_h,\pi_h) - \dot{b}(q,\tau q)(\du_h,\pi_h)
			\quad&\forall\ \pi_h \in P_h,\\
			\end{aligned}\right.
			\label{eq:d2Sh}
			\end{equation}
			with $(\tu_h,\tau p_h) = S_h'(q)(\tau q)$.
	\end{enumerate}

Like in the previous section, in order to study the convergence of the discrete quantities to their continuous counterparts, we introduce projection operators onto the discrete spaces. Since there will be no room for misunderstanding, to avoid redundant notation, all of them will be indicated by the same symbol $\Pi_h^r$, never minding if returning functions in $V_h,P_h,$ or $V_h\times P_h$.
\\
Referring to $X_h^r$, the following interpolation estimate is known (see e.g. \cite{Quarteroni2008a}, section 3.4.2), for $r\geq1,\ m=0,1$:
\beq
	|\varphi-\Pi_h^r\varphi|_{H^m(\O)}\leq ch^{r+1-m}|\varphi|_{H^{r+1}(\O)}.
\label{eq:interpX}
\eeq
The particular choice of $\mathds{P}_2-\mathds{P}_1$ couple in the spaces defined in \eqref{eq:spacesh}, leads us to assume the following regularity for the state variables and their shape derivatives:
\begin{ssmptn}
	For any $q\in\Q,\ \dq\in\delta Q$, any of $S(q),S'(q)(\dq),S''(q)(\dq,\dq)$ belong to $[H^3(\O)]^2\times H^2(\O)$ and the following inequalities hold:
	\beq
		\begin{aligned}
		\|S(q)\|_{[H^3(\O)]^2\times H^2(\O)}=&&\|\u\|_{[H^3(\O)]^2}&+\|p\|_{H^2(\O)}&&\leq c_1,\\
		\|S'(q)(\dq)\|_{[H^3(\O)]^2\times H^2(\O)}=&&\|\du\|_{[H^3(\O)]^2}&+\|\delta p\|_{H^2(\O)}&&\leq c_2\|\dq\|_{H^3(I)},\\
		\|S''(q)(\dq,\dq)\|_{[H^3(\O)]^2\times H^2(\O)}=&&\|\delta\du\|_{[H^3(\O)]^2}&+\|\delta\delta p\|_{H^2(\O)}&&\leq c_3\|\dq\|_{H^3(I)}^2.
		\end{aligned}
	\eeq
\label{ass:regolarita}
\end{ssmptn}
\begin{rmrk}
	In Appendix \ref{app:regolarita} we prove (Theorem \ref{th:regolaritaS}) the validity of Assumption \ref{ass:regolarita}, that involves suitable regularity assumptions on data of Stokes problem.
\end{rmrk}
Assumption \ref{ass:regolarita}, together with \eqref{eq:interpX}, yields the following estimate:
\beq
	\begin{aligned}
	\|\u-\Pi_h^2\u\|_V&+\|p-\Pi_h^1p\|_P &&\leq c_1h^2,\\
	\|\du-\Pi_h^2\du\|_V&+\|\delta p-\Pi_h^1\delta p\|_P &&\leq c_2h^2\|\dq\|_{H^3(I)},\\
	\|\delta\du-\Pi_h^2\delta\du\|_V&+\|\delta\delta p-\Pi_h^1\delta\delta p\|_P &&\leq c_3h^2\|\dq\|_{H^3(I)}^2,
	\end{aligned}
\label{eq:interp}
\eeq
which are crucial to obtain the following convergence result.
\begin{rmrk}
	Under regularity Assumption \ref{ass:regolarita}, one can afford the optimal convergence rate for $\mathds{P}_2-\mathds{P}_1$ discretization: lower regularity of the state variables would lead to a lower order on $h$ in \eqref{eq:interp}.
\end{rmrk}

The interpolation error estimates are once again the basis upon which we build our convergence result, which reads as follows:
\begin{lmm}
	For any $q_\sigma\in \Q_\sigma,\dq\in \delta Q$, the following convergence estimates hold:
	\begin{enumerate}[\quad(a)]
		\item $\|S(q_\sigma)-S_h(q_\sigma)\|_{V\times P} \leq ch^2,$
		\item $\|S'(q_\sigma)(\dq)-S'_h(q_\sigma)(\dq)\|_{V\times P} \leq ch^2\|\dq\|_{H^3(I)},$
		\item $\|S''(q_\sigma)(\dq,\dq)-S''_h(q_\sigma)(\dq,\dq)\|_{V\times P} \leq ch^2\|\dq\|_{H^3(I)}^2$.
	\end{enumerate}
\label{th:SShcont}
\end{lmm}
\begin{proof} Since the discrete problems \eqref{eq:StokesdebTh}-\eqref{eq:d2Sh} fulfill the same properties as the continuous ones, we have that Theorem \ref{th:dS} on the boundedness of the continuous solution operator $S$ is true also for the discrete operator $S_h$ and its derivatives. Hinging upon this result and Assumption \ref{ass:regolarita}, we fix some $q_\sigma\in\Q_\sigma,\dq\in\delta Q$ and proceed according to the following steps.

We first prove (\emph{a}). From \cite{Girault1986} and the independence of the continuity, coercivity and LBB constants from $q_\sigma$ and $h$, we can obtain the classical convergence result for a saddle-point problem, i.e.,
	\beq
		\|S(q_\sigma)-S_h(q_\sigma)\|_{V\times P} \leq c(\|\u-\Pi^2_h\u\|_V+\|p-\Pi^1_hp\|_P) \leq ch^2,
	\eeq
with the last inequalities exploiting interpolation error estimate \eqref{eq:interp}. 

We now proceed to prove (\emph{b}). We set $(\u,p)=S(q)$, $(\du,\delta p)=S'(q)(\delta q)$, $(\delta\du,\delta\delta p)=S''(q)(\dq,\dq)$, with subscript $\cdot_h$ denoting the correspondent discrete quantities, and we introduce the ``intermediate derivative'' $(\delta\hat{\u}_h,\delta\hat{p}_h)$, solution in $V_h\times P_h$ of the following problem:
\footnote{The problem here introduced is a combination of problem \eqref{eq:dS} for $S'(q_\sigma)(\dq)$ and its discrete counterpart \eqref{eq:dSh}: we solve a discrete problem in spaces $V_h,P_h$, with the first equation being the same as in \eqref{eq:dS}, and the second one as in \eqref{eq:dSh}.}
	\begin{equation}
		\left\{
		\begin{aligned}
		a(q)(\delta\hat{\u}_h,\mathbf{v}_h) + b(q)(\mathbf{v}_h,\delta\hat{p}_h) &= \dot{F}(q,\dq)(\mathbf{v}_h)+&\\
		& - \dot{a}(q,\dq)(\u,\mathbf v_h) - \dot{b}(q,\dq)(\mathbf v_h, p)\qquad&\forall\ \mathbf{v}_h\in  V_h,\\
		b(q)(\delta\hat{\u}_h,\pi_h) &= \dot{G}(q,\dq)(\pi_h)-\dot{b}(q,\dq)(\u_h,\pi_h)\qquad&\forall\ \pi_h \in P_h.\\
		\end{aligned}\right.\\
	\label{eq:new_dSinterm}
	\end{equation}
	Thanks to \eqref{eq:new_dSinterm}, we can separate the error due to the discretization of the problem on $S'(q)(\dq)$ from the one that is inherited from the discretization of $S(q)$. Using triangular inequality yields
	\beq
	\begin{split}
		 \|S'(&q)(\dq)-S_h'(q)(\dq)\|_{V\times P} \leq \\
		&\leq\|S'(q)(\dq)-(\delta\hat{\u}_h,\delta\hat{p}_h)\|_{V\times P} + \|(\delta\hat{\u}_h,\delta\hat{p}_h)-S_h'(q)(\dq)\|_{V\times P} = \\
		 &= \|\du-\delta\hat{\u}_h\|_V + \|\delta\hat{\u}_h-\du_h\|_V + \|\delta p - \delta\hat{p}_h\|_P + \|\delta\hat{p}_h-\delta p_h\|_P.
	\end{split}
	\label{eq:new_SSh1}
	\eeq
	Considering the first term in \eqref{eq:new_SSh1}, we have that, for any $\mathbf w_h\in V_h$,
	\beq\begin{split}
		&\alpha_c\|\nabla\du - \nabla\delta\hat{\u}_h\|^2 \leq a(q_\sigma)(\du - \delta\hat{\u}_h,\du - \delta\hat{\u}_h) =\\
		&=a(q_\sigma)(\du - \delta\hat{\u}_h,\du-\mathbf w_h) +a(q_\sigma)(\du-\delta\hat{\u}_h,\mathbf w_h-\delta\hat{\u}_h)=\\
		&=a(q_\sigma)(\du - \delta\hat{\u}_h,\du-\mathbf w_h) -b(q_\sigma)(\mathbf w_h-\delta\hat{\u}_h,\delta p-\delta\hat{p}_h),
	\end{split}\label{eq:new_du-duinterminit}
	\eeq
	with the equality holding thanks to the fact that the first equations in \eqref{eq:dS} and \eqref{eq:new_dSinterm} share the same right-hand side. 
Since \eqref{eq:new_du-duinterminit} holds for every $\mathbf w_h \in V_h$, it still holds if we take the infimum w.r.t.\ $\mathbf w_h$. For the first term of the right-hand side we get
	\beq
		\begin{split}
			\inf_{\mathbf w_h\in V_h}a(q_\sigma&)(\du - \delta\hat{\u}_h,\du-\mathbf w_h)\leq
a(q_\sigma)(\du - \delta\hat{\u}_h,\du-\Pi_h\du)\leq \\
			&\leq M\|\nabla\du-\nabla\delta\hat{\u}_h\|\|\nabla\du-\nabla\,\Pi_h\du\|\leq ch^2M\|\nabla\du-\nabla\delta\hat{\u}_h\|\|\dq\|_{H^3(I)},
		\end{split}
	\label{eq:new_aqdu-duinterm}
	\eeq
	where we employed the interpolation error estimate \eqref{eq:interp} and the boundedness of $\|\du\|_{H^2(I)}$ (due to Assumption \ref{ass:regolarita}). 
	Instead, taking $\mathbf w_h=\delta\hat{\u}_h$ in the second term yields
	\beq
		\inf_{\mathbf w_h\in V_h}[-b(q_\sigma)(\mathbf w_h-\delta\hat{\u}_h,\delta p-\delta\hat{p}_h)]\leq0.
	\label{eq:new_bqdu-duinterm}
	\eeq
	Using \eqref{eq:new_aqdu-duinterm} and \eqref{eq:new_bqdu-duinterm} in \eqref{eq:new_du-duinterm} and dividing both sides by $\alpha_c\|\nabla\du - \nabla\delta\hat{\u}_h\|$, we eventually obtain
	\beq
		\|\nabla\du - \nabla\delta\hat{\u}_h\|\leq c\frac{M}{\alpha_c}h^2\|\dq\|_{H^3(I)}.
	\label{eq:new_du-duinterm}
	\eeq
	The second term in \eqref{eq:new_SSh1} can be estimated using the problems \eqref{eq:StokesdebTh} and \eqref{eq:new_dSinterm}, fulfilled by $\du_h, \delta\widehat{\u}_h$, together with the coercivity of $a$ and the continuity of the forms involved in such problems. We can thus obtain:
	\beq\begin{split}
		&\alpha_c\|\nabla\delta\hat{\u}_h-\nabla\du_h\|^2 \leq a(q_\sigma)(\delta\hat{\u}_h-\du_h,\delta\hat{\u}_h-\du_h) =\\
		&= -\dot{a}(q_\sigma,\dq)(\u-\u_h,\delta\hat{\u}_h-\du_h)-\dot{b}(q_\sigma,\dq)(\delta\hat{\u}_h-\du_h,p-p_h) +\\
			&\qquad - b(q_\sigma)(\delta\hat{\u}_h-\du_h,\delta\hat{p}_h-\delta p_h)\leq\\
		&\leq c\|\dq\|_{H^2(I)}(\|\nabla\u-\nabla\u_h\|+\|p-p_h\|)\|\nabla\delta\hat{\u}_h-\nabla\du_h\|, 
	\end{split}\label{eq:new_aqduinterm-duh}
	\eeq
	where the last inequality holds because $b(q_\sigma)(\delta\hat{\u}_h,\pi_h)=b(q_\sigma)(\du_h,\pi_h)\quad\forall\pi_h\in P_h$.
	After dividing by $\|\nabla\delta\hat{\u}_h-\nabla\du_h\|$ both sides of \eqref{eq:new_aqduinterm-duh}, the right-hand side can be controlled as in the first point of the present Lemma, leading to
	\beq
		\|\nabla\delta\hat{\u}_h-\nabla\du_h\|\leq ch^2\|\dq\|_{H^3(I)}.
	\label{eq:new_duinterm-duh}
	\eeq

	Now we have to deal with pressure error terms in \eqref{eq:new_SSh1}: taking a generic $\pi_h\in P_h$, the first term can be split as follows:
	\beq
		\|\delta p -\delta p_h\| \leq \|\delta p - \pi_h\| + \|\pi_h-\delta\hat{p}_h\|.
	\label{eq:new_dp-dphinit}
	\eeq
	We remark that, since inequality \eqref{eq:new_dp-dphinit} holds for any $\pi_h\in P_h$, it holds also taking the infimum w.r.t.\ $\pi_h$. The infimum of the first term is directly controlled by $ch^2\|\dq\|_{H^2(I)}$ thanks to the interpolation error estimate \eqref{eq:interp} and the boundedness of $\|\delta p\|_{H^2(\O)}$ asserted in Theorem \ref{th:regolaritaS}. The second term goes to zero when passing to the infimum, since $\delta p_h\in P_h$.

	Finally, for the last term in \eqref{eq:new_SSh1} we exploit LBB condition \eqref{eq:LBB} and proceed as follows:
	\beq\begin{split}
		&\|\delta\hat{p}_h-\delta p_h\|\leq\sup_{\v_h\in V_h}\frac{b(q_\sigma)(\v_h,\delta\hat{p}_h-\delta p_h)}{\hat{\beta}\|\nabla\v_h\|}=\\
		&=\sup_{\v_h\in V_h}\frac{-\dot{a}(q_\sigma,\dq)(\u-\u_h,\v_h)-\dot{b}(q_\sigma,\dq)(\v_h,p-p_h) - a(q_\sigma)(\delta\hat{\u}_h-\du_h,\v_h)}{\hat{\beta}\|\nabla\v_h\|}\leq\\
		&\leq\frac{1}{\hat{\beta}}\left[c\|\dq\|_{H^2(I)}\left(\|\nabla\u-\nabla\u_h\|+\|p-p_h\|\right)+ M\|\nabla\delta\hat{\u}_h-\nabla\du_h\|\right].
	\end{split}\label{eq:new_dpinterm-dph}
	\eeq
	From estimate \eqref{eq:new_duinterm-duh} and point (\emph{a}) of the present lemma, we get the desired bound, i.e. $ch^2\|\dq\|_{H^3(I)}$.\\
	Collecting the estimates for the four terms in \eqref{eq:new_SSh1} yields the validity of point (\emph{b}).

	Finally, we prove (\emph{c}), employing the regularity result for $S''(q)(\dq,\dq)$ given at the third point of Assumption \ref{ass:regolarita}. The only difference from the previous point is the ``intermediate derivative'' \linebreak[4]$(\delta\delta\check{\u}_h,\delta\delta\check{p}_h)\in V_h\times P_h$, defined as the solution of the following problem:
	\begin{equation}
			\left\{
			\begin{aligned}
			&\begin{aligned}
			a(q)&(\delta\delta\check{\u}_h,\mathbf{v}_h) + b(q)(\mathbf{v}_h,\delta\delta\check{p}_h) = \\
			&=\ddot{F}(q,\dq,\dq)(\mathbf{v}_h) - \ddot{a}(q,\dq,\dq)(\u,\mathbf v_h) - \ddot{b}(q,\dq,\dq)(\mathbf v_h, p)+\\
			&- 2\,\dot{a}(q,\dq)(\du,\mathbf v_h) - 2\,\dot{b}(q,\dq)(\mathbf v_h, \delta p)
			\end{aligned}\quad&\forall\ \mathbf{v}_h\in  V_h,\\
			&\begin{aligned}
			b(q)&(\delta\delta\check{\u}_h,\pi_h) = \ddot{G}(q,\dq,\dq)(\pi_h)-\ddot{b}(q,\dq,\dq)(\u_h,\pi_h) +\\
			&- 2\,\dot{b}(q,\dq)(\du_h,\pi_h) 
			\end{aligned}\quad&\forall\ \pi_h \in P_h.\\
			\end{aligned}\right.
	\label{eq:new_d2Sinterm}
	\end{equation}
	All the previous steps performed to estimate $S'-S_h'$ can be easily adapted to the present context.

\end{proof}

A direct consequence of the previous lemma is the following convergence result for the discrete functional.
\begin{lmm}
	$\forall q_\sigma\in\Q_\sigma,\dq\in \delta Q$ it holds
	\begin{enumerate}[\quad(a)]
	\item $|j(q_\sigma)-j_h(q_\sigma)|\leq ch^2,$\label{punto1}
	\item $|j'(\qs)(\dq)-j_h'(\qs)(\dq)|\leq ch^2\|\dq\|_{H^3(I)},$
	\item $|j''(\qs)(\dq,\dq)-j_h''(\qs)(\dq,\dq)|\leq ch^2\|\dq\|_{H^3(I)}^2.$
	\end{enumerate}
\label{th:j-jh}
\end{lmm}
\begin{proof}
	Let us fix a  $\qs\in\Q_\sigma,\ \dq\in \delta Q$ and define $(\u,p)=S(\qs),\ (\du,\delta p)=S'(\qs)(\dq),$ $(\delta\du,\delta\delta p)=S''(\qs)(\dq,\dq)$.

	Let us first prove (\emph{a}). It it easy to show that the following holds
	\beq\begin{split}
		&|j(\qs)-j_h(\qs)| = |((\nabla\u-\nabla\u_h)\,A_q,\nabla\u+\nabla\u_h)|\leq\\
		&\leq\|A_q\|_\infty(\|\nabla\u\|+\|\nabla\u_h\|)\|\nabla\u-\nabla\u_h\|\leq ch^2,
	\end{split}\eeq
	where the last inequality employs the boundedness of $A_q,\nabla\u,\nabla\u_h$ and Lemma \ref{th:SShcont}.

	Now we prove (\emph{b}), according to the following steps:
	\beq\begin{split}
	|j'(\qs)(\dq)-j_h'(\qs)(\dq)|&\leq|((\nabla\u-\nabla\u_h)\,A'_{q,\dq},\nabla\u+\nabla\u_h)|+\\
	&\quad+2|((\nabla\du-\nabla\du_h)\,A_q,\nabla\u)| + 2|(\nabla\du_h\,A_q,\nabla\u-\nabla\u_h)|\leq\\
	&\leq ch^2\|\dq\|_{H^3(I)}.
	\end{split}\eeq
	Indeed, it holds $\|A'_{q,\dq}\|_\infty\leq c\|\dq\|_{H^2(I)}$ (see Proposition \ref{th:coeff}) while $\|\nabla\u\|$ and $\|\nabla\du_h\|$ are controlled thanks to the continuous and discrete versions of Theorem \ref{th:dS}, and the discretization error terms are bounded through Lemma \ref{th:SShcont}.

	Finally, we prove (\emph{c}), as follows:
	\beq\begin{split}
	&|j''(\qs)(\dq,\dq)-j_h''(\qs)(\dq,\dq)|\leq|((\nabla\u-\nabla\u_h)\,A''_{q,\dq,\dq},\nabla\u+\nabla\u_h)|+\\
	&+4|((\nabla\du-\nabla\du_h)\,A'_{q,\dq},\nabla\u)| + 4|(\nabla\du_h\,A'_{q,\dq},\nabla\u-\nabla\u_h)|+\\
	&+2|((\nabla\du-\nabla\du_h)\,A_q,\nabla\du+\nabla\du_h)|+\\
	&+2|((\nabla\delta\du-\nabla\delta\du_h)\,A_q,\nabla\u)| + 2|(\nabla\delta\du_h\,A_q,\nabla\u-\nabla\u_h)|.
	\end{split}\eeq
	To bound the terms not involving $\delta\du$ and $\delta\du_h$, one can employ Proposition \ref{th:coeff} to handle the matrix terms, together with similar techniques already used to prove (\emph{a}) and (\emph{b}). To bound the last two terms, we have to apply Lemma \ref{th:SShcont}, point \emph{c}, and Theorem \ref{th:regolaritaS} in order to provide estimates for \linebreak[4]$\|\nabla\delta\du-\nabla\delta\du_h\|$ and $\|\nabla\delta\du_h\|$.
\end{proof}

\ \\

Finally, collecting the previous results, we can prove the main result of this section.
\begin{thrm}[A priori convergence estimates]
	Let Assumptions \ref{ass:noeps}, \ref{ass:jpos} and \ref{ass:regolarita} hold. Then, denoted by $\overline{q}$ a local solution of \eqref{eq:minj}, there exists a sequence $\{\overline{q}_{\sigma,h}\}_{\sigma,h>0}$ of local optimal solution of the discrete problem
	\beq
		\min_{\qs\in\Q_\sigma} j_h(\qs),
	\label{eq:minJsigmah}
	\eeq
such that
	\begin{gather}
	\|\overline{q}-\overline{q}_{\sigma,h}\|_{H^3(I)}=\mathcal{O}(\sigma^2+h^2),\\
	\|S(\overline{q})-S_h(\overline{q}_{\sigma,h})\|_{V\times P}=\mathcal{O}(\sigma^2+h^2),\\
	|j(\overline{q})-j_h(\overline{q}_{\sigma,h})|=\mathcal{O}(\sigma^2+h^2).
	\end{gather}
\label{th:conv}
\end{thrm}
\begin{proof}
	Let $\overline{q}_\sigma,\overline{q}_{\sigma,h}$ denote the optimal controls for the semi-discrete problem \eqref{eq:minJsigma} and the fully discretized problem \eqref{eq:minJsigmah}, respectively.
	The Mean Value Theorem ensures the existence of $t\in(0,1)$ such that, with $\xi=t\overline{q}_\sigma+(1-t)\overline{q}_{\sigma,h}$, we have
	\beq
		j_h'(\overline{q}_\sigma)(\dq_\sigma)-j_h'(\overline{q}_{\sigma,h})(\dq_\sigma) = j''_h(\xi)(\dq_\sigma,\overline{q}_\sigma-\overline{q}_{\sigma,h}).
	\label{eq:jxi}
	\eeq
	Applying Lemma \ref{th:314} and taking $\overline{q}_\sigma-\overline{q}_{\sigma,h}$ as a variation, we get:
	\beq
		\begin{split}
		\frac{\delta_2}{2}&\|\overline{q}_\sigma-\overline{q}_{\sigma,h}\|_{H^3(I)}^2\leq j''(\xi)(\overline{q}_\sigma-\overline{q}_{\sigma,h},\overline{q}_\sigma-\overline{q}_{\sigma,h})\leq\\
		&\leq j_h''(\xi)(\overline{q}_\sigma-\overline{q}_{\sigma,h},\overline{q}_\sigma-\overline{q}_{\sigma,h})+\\
		&\qquad\qquad+ | j''(\xi)(\overline{q}_\sigma-\overline{q}_{\sigma,h},\overline{q}_\sigma-\overline{q}_{\sigma,h})- j_h''(\xi)(\overline{q}_\sigma-\overline{q}_{\sigma,h},\overline{q}_\sigma-\overline{q}_{\sigma,h})|\leq\\
		&\leq j_h'(\overline{q}_\sigma)(\overline{q}_\sigma-\overline{q}_{\sigma,h})- j_h'(\overline{q}_{\sigma,h})(\overline{q}_\sigma-\overline{q}_{\sigma,h}) + c_1h^{2}\|\overline{q}_\sigma-\overline{q}_{\sigma,h}\|_{H^3(I)}^2,
		\end{split}
	\label{eq:passaggiqsigma-qsigmah}
	\eeq
	where the last inequality is obtained by \eqref{eq:jxi} and Lemma \ref{th:j-jh}(\emph{c}). Using the fact
that $j'_h(\overline{q}_{\sigma,h})(\overline{q}_\sigma-\overline{q}_{\sigma,h})=j'(\overline{q}_\sigma)(\overline{q}_\sigma-\overline{q}_{\sigma,h})=0$ in the right-hand side of \eqref{eq:passaggiqsigma-qsigmah} and then applying Lemma \ref{th:j-jh}(\emph{b}), we obtain:
	\beq
		\begin{split}
		\frac{\delta_2}{2}&\|\overline{q}_\sigma-\overline{q}_{\sigma,h}\|_{H^2(I)}^2\leq j_h'(\overline{q}_\sigma)(\overline{q}_\sigma-\overline{q}_{\sigma,h})- j'(\overline{q}_\sigma)(\overline{q}_\sigma-\overline{q}_{\sigma,h}) + c_1h^{2}\|\overline{q}_\sigma-\overline{q}_{\sigma,h}\|_{H^3(I)}^2\leq\\
		&\leq c_2h^2\|\overline{q}_\sigma-\overline{q}_{\sigma,h}\|_{H^3(I)} + c_1h^{2}\|\overline{q}_\sigma-\overline{q}_{\sigma,h}\|_{H^3(I)}^2.
		\end{split}
	\eeq
	Therefore, for sufficiently small $h$, i.e. for
	\beq
		h\leq\left(\frac{\delta_2}{2c_1}\right)^{1/2},
	\eeq
	the following convergence error estimate holds:
	\beq
		\|\overline{q}-\overline{q}_{\sigma,h}\|_{H^3(I)}\leq\|\overline{q}-\overline{q}_\sigma\|_{H^3(I)}+\|\overline{q}_\sigma-\overline{q}_{\sigma,h}\|_{H^3(I)}=\mathcal{O}(\sigma^2+h^2).
	\eeq
	This result yields the second point of thesis, since
	\beq
		\|S(\overline{q})-S_h(\overline{q}_{\sigma,h})\|_{V\times P}\leq\|S(\overline{q})-S(\overline{q}_{\sigma,h})\|_{V\times P}+\|S(\overline{q}_{\sigma,h})-S_h(\overline{q}_{\sigma,h})\|_{V\times P},
	\label{eq:S-Shsplit}
	\eeq
	and the desired estimate for $S_h$ follows from applying Lemmas \ref{th:SLip} and \ref{th:SShcont} to the two terms at right-hand side of \eqref{eq:S-Shsplit}. An analogous argument, using Lemmas \ref{th:jLip} and \ref{th:j-jh}, yields the estimate for $j_h$.
\end{proof}

\section{Numerical results}\label{sec:algores}
In this section, we present two sets of numerical results. 
The numerical implementation has been carried out basing on the \fenics\ project (see \cite{Logg} and \url{http://fenicsproject.org}), and the optimal solution is obtained iteratively, using %
the following gradient method \cite{Nocedal1999}:

\mybox{
	Given $q_{old}$ from the previous iteration,\\
	set the descent step length $\varepsilon$ to the initial value $\widehat{\varepsilon}>0$. Then, \hfill
	\begin{enumerate}
		\item solve state and adjoint problems in order to obtain $(\u,p),(\z,s)$
		\item build $\nabla j(q_{old})$
		\item project $\nabla j(q_{old})$ on the set of admissible variations, obtaining $G$
		\item restrict $G$ on $\Gamma_0$ and then map it to $I$, to get $g$
		\item back-tracking: set $q_{new} = q_{old}-\widehat{\varepsilon}g$\\
		 	while $j(q_{new}) > j(q_{old})$ and $\varepsilon>\varepsilon_{min}$ do: 
		\begin{enumerate}
				\item update: $q_{new} = q_{old} - \varepsilon g$
				\item $\varepsilon = \varepsilon/2$
		\end{enumerate}
	\end{enumerate}
\label{iter_grad}
}{Gradient method iteration}
In general, the functional gradient $\nabla j(q_{old})$, obtained as in Lemma \ref{th:Hadamard}, is not an admissible variation, since one cannot prove the existence of some $\varepsilon>0$ such that ${q=q_{old}-\varepsilon\nabla j(q_{old})}$ satisfies
\beq
	q(0)=q(1)=0.
\eeq
This is why in the gradient method the projection step (3) is required. The gradient $\nabla j(q_{old})$ is projected onto $H^1_{\partial \O\setminus\Gamma_0}(\O)$ solving the following problem:
\beq
\label{eqN:G}
	\begin{cases}
		-\Delta G + G = 0,\qquad\ &\text{in } \O,\\
		G = 0,\qquad\ &\text{on } \partial\O\setminus\Gamma_0,\\
		- \partial_{\mathbf n} G = -\nabla j(q_{old}),\qquad\ &\text{on } \Gamma_0.
	\end{cases}
\eeq
Then, step (4) of the algorithm reduces $G$, defined on $\O$, to a function $g$ belonging to the space of controls.

The results obtained by the application of the above algorithm to the shape optimization problem \eqref{eq:minj} are now presented and discussed. Two different functionals will be considered in the two test cases.
\begin{rmrk}
	We remark that we use finite element discretization, with $\mathds{P}_2-\mathds{P}_1$ pair for state velocity and pressure and with piecewise linear basis functions for the control. As we will see, even if the polynomial degree for controls is not as high as assumed in the derivation of a priori estimates, the numerical results comply the theoretical ones. In these numerical tests, we consider a unique discretization parameter, i.e. we set $\sigma=h$.
\end{rmrk}

\subsection{Test case 1}

In this first test case, we take into account the following functional:
\beq
	\widetilde{j}(q) = \int_\Oq |\nabla\widetilde{\u}|^2\,d\Omega + \alpha\int_{\Gamma_q}d\Gamma + \beta\left(\int_Iq(x)\,dx \ -\ \overline{V}\right)^2.
\eeq
Its counterpart on the pulled-back formulation \eqref{eq:StokesdebT} reads
\beq
	j(q) = \int_\O |\nabla\u DT_q^{-1}|^2\,d\Omega + \alpha\int_I\sqrt{1+(q'(x))^2}\,dx + \beta\left(\int_Iq(x)\,dx \ -\ \overline{V}\right)^2.
\label{eqN:j}
\eeq
The gradient of this functional is given by
\beq
\begin{split}
	\nabla j(q) &=[\nabla{\u}DT_q^{-1}DT_q^{-T}\mathbf n]\cdot\left[(\nu\nabla\z -{\u})DT_q^{-1}DT_q^{-T}\mathbf n\right] |_{\Gamma_0}+\\
	&- 2\alpha\frac{q''}{1+(q')^2}+2\beta\left(\int_Iq(x)dx - \overline{V}\right).
\end{split}
\label{eqN:gradj}
\eeq
The regularization term considered in \eqref{eqN:j} 
is often used in literature (see, e.g., \cite{Dogan2007,Morin2012}) and it consists in the penalization of the perimeter of the moving portion $\Gamma_q$ of the domain boundary.
This new term 
is simpler to handle than the curvature term $\|q''\|_{L^2(I)}^2$: indeed, using the original term would require the introduction of a further adjoint problem, to extract the Riesz representative in $L^2(I)$ of $\dq\mapsto(q'',\dq'')_I$. Moreover, the perimeter term can be supposed to generally act in the same way as the curvature term, since a shorter perimeter corresponds to less oscillations, and vice versa.

\begin{figure}[h!]
	\subfloat[Initial configurations\label{fig:initall}]{\includegraphics[width=0.45\textwidth]{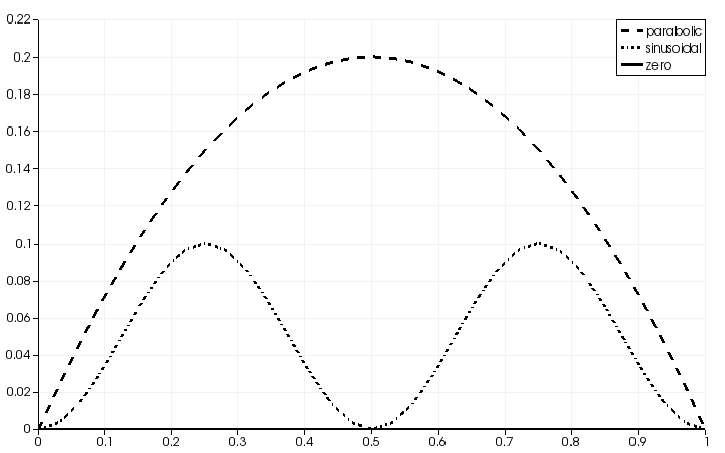}}\qquad
	\subfloat[Optimal controls\label{fig:endall}]{\includegraphics[width=0.45\textwidth]{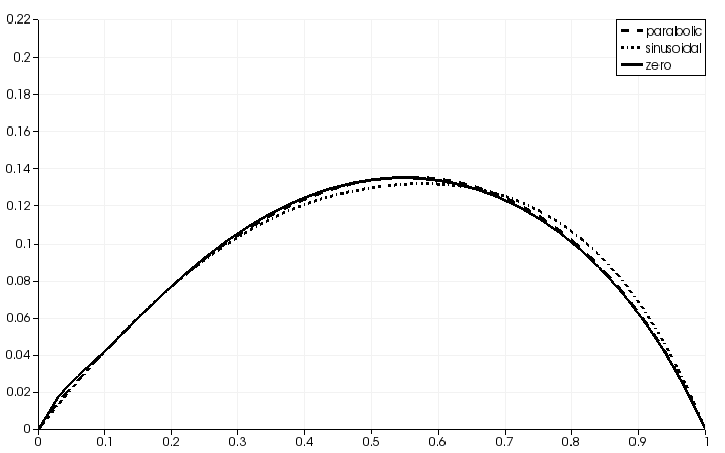}}\\
\caption{Independence of the optimal control from the initial configuration, for $\alpha=10$, $\beta=10\,000$, $\overline{V}=[0.7$ times the initial area of the parabolic case$]$.}
\label{fig:all}
\end{figure}

We first analyze the dependence of the optimal solution on the initial configuration. We considered three different initial solutions, defined by a parabolic function ($q(x)=0.2[1-4(x-0.5)^2]$), a sinusoidal function ($q(x)=0.1\sin(2\pi x)^2$), and the flat function ($q(x)\equiv0$).
As shown in \figref{fig:all}, the optimal control obtained are very close, starting from different initial controls. The final configurations in \figref{fig:endall} are reached in less than 10 iterations, with $\hat{\varepsilon}=0.1, \varepsilon_{min}=10^{-8},$ and the reaching of $\varepsilon \leq \varepsilon_{min}$ as the stop criterion on the iterations of the gradient method.

\begin{figure}[h]
	\subfloat[Varying $\alpha$; $\beta=0$\label{fig:endParabPer}]{\includegraphics[width=0.45\textwidth]{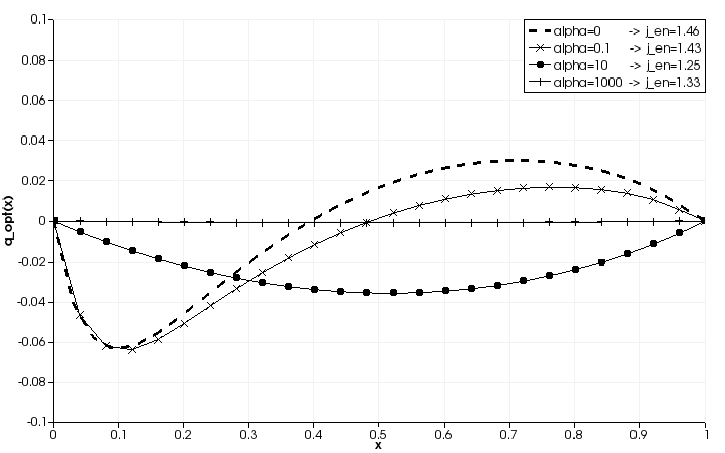}}\qquad
	\subfloat[Varying $\beta$; $\alpha=0$\label{fig:endParabVolq}]{\includegraphics[width=0.45\textwidth]{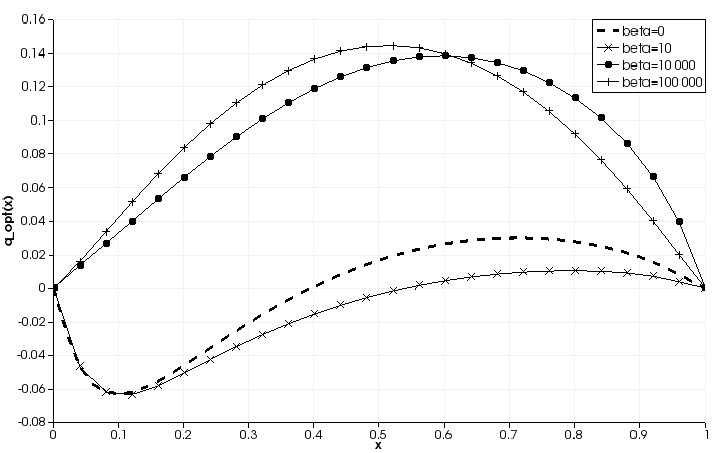}}
\caption{Final controls obtained by the optimization algorithm for different values of the penalty parameters ($j_{en}$ is the energetic term of the functional $j$)}
\label{fig:penalty}
\end{figure}

The dependence of the solution on the value of the penalty parameters has also been analyzed, 
starting from the parabolic configuration in \figref{fig:initall}.
Concerning parameter $\alpha$, a minimum value has to be exceeded in order to prevent the gradient method from converging to a local, sub-optimal minimum. Indeed, \figref{fig:endParabPer} shows that for lower values of $\alpha$, oscillating controls are found at the end of the optimization algorithm, though the value of the functional in such configurations is higher than the ones corresponding to $\alpha=10,1000$. Moreover, a maximum value must not be exceeded, otherwise the regularization parameter dominates too much in the total functional value, leading to a nearly flat optimal control. About parameter $\beta$, instead, we only need it to be greater than a minimum threshold, in order to sufficiently express the volume constraint. Under these considerations, \figref{fig:penalty} shows that the values $\alpha=10, \beta=10\,000$, considered in the previous test, are suitable for a proper expression of the two penalty terms.

\FloatBarrier
\subsection{Test case 2}

In this section, we report a numerical convergence analysis, carried out to validate the a priori error estimates proved in Theorem \ref{th:conv}.
For this purpose, we would like to have an exact solution as a reference point. To this end, we take into account the following functional:
\beq
	\widetilde{j}(q) = \int_\Oq |\nabla\widetilde{\u}-\nabla\widetilde{\u}_d|^2\,d\Omega + \frac{\alpha}{2}\int_{\Gamma_q}d\Gamma,
\eeq
with its pulled-back counterpart given by
\beq
	j(q) = \int_\O(\nabla\u-\nabla\u_d)A_q(\nabla\u-\nabla\u_d)d\Omega+ \frac{\alpha}{2}\int_I\sqrt{1+(q'(x))^2}\,dx.
\label{eqN:jnew}
\eeq

The velocity $\widetilde{\u}_d$ is obtained solving the Stokes problem on a domain $\Omega_{q_d}$, identified by the given control function
\beq
	q_d = 0.1+0.1\cos(2\pi (x-0.5)),
\eeq
and $\u_d=\widetilde{\u}_d\circ T_{q_d}^{-1}$.

Indeed, if no penalty terms are active, the minimum for this functional is zero, and it is reached for $q=q_d$. The functional \eqref{eqN:jnew} is a slight generalization of the functional \eqref{eq:jtildedef}, and the theoretical results presented in the previous sections can be easily generalized to the new functional.


%
Following the steps of Section \ref{sec:ottimalita}, we can derive an expression for the shape gradient in $q$:
\beq
\begin{split}
	\nabla j(q) &=-2\alpha\frac{q''(x)}{1+(q'(x))^2}+\\
	&+[(\nabla{\u}-\nabla\u_d)DT_q^{-1}DT_q^{-T}\mathbf n]\cdot\left[(\nu\nabla\z_{\u_d} -{\u}+\u_d)DT_q^{-1}DT_q^{-T}\mathbf n\right],
\end{split}
\label{eqN:gradjnew}
\eeq
where $\z_{\u_d}$ is the adjoint velocity variable, solution of a problem obtained from a minimal modification of \eqref{eq:PAHad}, replacing any occurrence of $\tilde{\u}$ with $\tilde{\u}-\tilde{\u}_d$.\\

Based on the functional defined in \eqref{eqN:jnew}, different spatial convergence tests have been carried out, taking four specific values for perimeter penalty coefficient $\alpha$, namely $\alpha=0,0.01,0.1,1$. 

\begin{figure}[h]
	\subfloat[$\alpha=0$]{\includegraphics[width=0.45\textwidth]{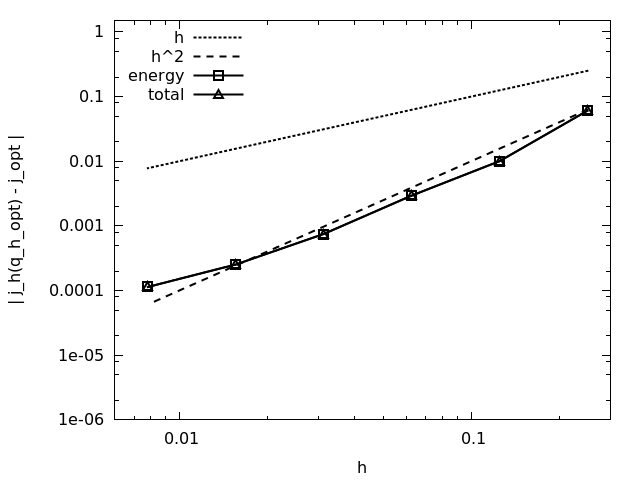}}\qquad
	\subfloat[$\alpha=0.01$]{\includegraphics[width=0.45\textwidth]{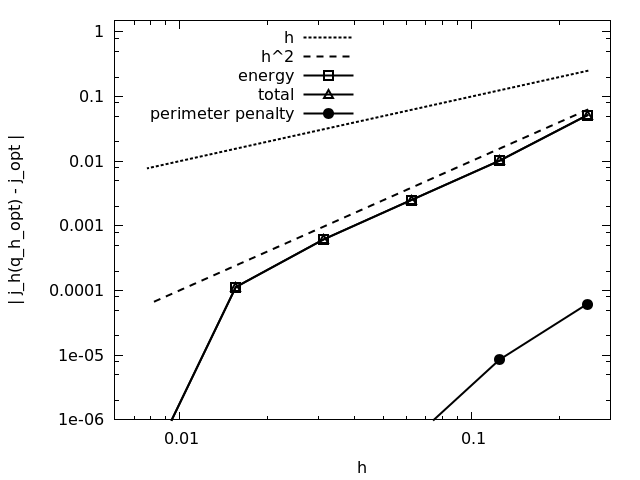}}\\
	\subfloat[$\alpha=0.1$]{\includegraphics[width=0.45\textwidth]{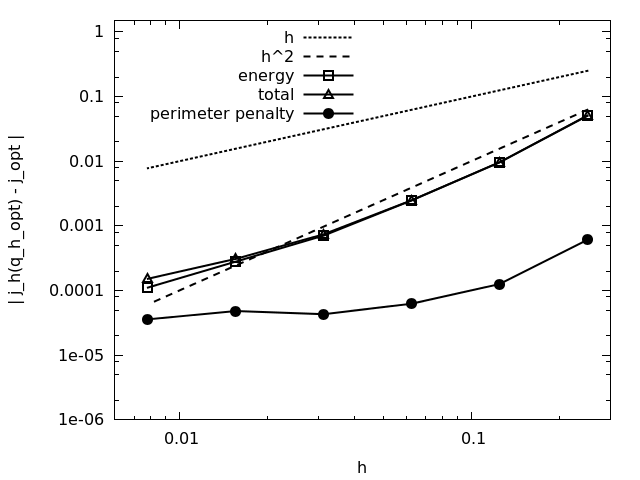}}\qquad
	\subfloat[$\alpha=1$]{\includegraphics[width=0.45\textwidth]{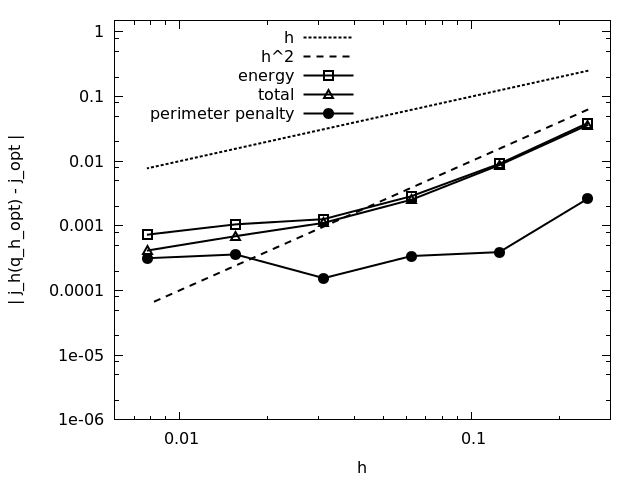}}
\caption[Spatial convergence of functional values]{Spatial convergence of discrete functional value $j_h(q_{h,opt})$ to its reference value $j(q_{opt})$. Each term of the functional is presented w.r.t.\ its corresponding term in $j(q_{opt})$, which is known for $\alpha=0$, and obtained by Richardson extrapolation for $\alpha\neq0$.}
\label{fig:convergence}
\end{figure}


The results reported in \figref{fig:convergence} are in agreement with the a priori estimates of the convergence error proved in Theorem \ref{th:conv}, since an approximately quadratic convergence order is obtained, for a broad spectrum of values of $h$. 
%
However, for $h\to0$, the graphs in \figref{fig:convergence} show a sort of saturation bending. A reason for this can be found in the stopping criterion of the optimization algorithm and in the lower bound imposed on the descent step length, that introduce a finite error. This influence is amplified as $\alpha$ grows, to the point of polluting the convergence behaviour, hence we do not report results for $\alpha>1$.



\section*{Conclusions}

In this paper, we have studied a shape optimization problem, namely the minimization of the total energy dissipation for the low-Reynolds flow of a viscous, incompressible fluid, modeled by two-dimensional, steady Stokes equations. After the definition of the problem and the admissible set of control functions, we have reformulated the problem onto a reference domain, by means of a control-dependent map. The well-posedness of the transformed problem has been inspected, and particular attention has been devoted to the inf-sup condition for the form $b(q)$, obtaining a control-independent lower bound for the inf-sup constant. The existence of an optimal solution has also been proved, for the minimization problem at hand, and corresponding first order optimality conditions have been provided.

After the inspection of some differentiability properties of the state solution operator, a FEM discretization of the problem has been introduced. For this discretization, a priori error estimates have been derived, showing a quadratic convergence rate. To our best knowledge, this is the first result about convergence rates obtained for the discretization of Stokes problem in a shape optimization environment. Numerical tests have been performed to assess the validity of the theoretical results.
%

\FloatBarrier
\appendix
\begin{appendices}
\section{Additional regularity}\label{app:regolarita}

In this Appendix we want to show a possible way to derive the regularity properties stated in Assumption \ref{ass:regolarita}, starting from suitable requests on data and a regularity result on Stokes problem with mixed boundary conditions.

At first, let us state a preliminary result about the transformation of norms defined on the reference domain ($\O$) and on the physical one ($\Oq$).
\begin{lmm}
	Let $k\in\mathds{N}$ be fixed, $\varphi\in H^{k}(\O)$ and $q\in W^{k,\infty}(I)$. It holds that
	$$\varphi\circ T_q^{-1}\in H^k(\Oq),\quad c_1\|q\|_{W^{k,\infty}(I)}\|\varphi\circ T_q^{-1}\|_{H^k(\Oq)}\leq\|\varphi\|_{H^k(\O)}\leq c_2\|q\|_{W^{k,\infty}(I)}\|\varphi\circ T_q^{-1}\|_{H^k(\Oq)}.$$
	Vice versa, it holds that $\widetilde{\varphi}\in H^k(\Oq)$ implies $\widetilde{\varphi}\circ T_q\in H^k(\O)$, together with similar inequalities.
\label{th:HkTq}
\end{lmm}
In connection with this lemma, we restrict a little the set of admissible controls.
	From now on, the definition of $\Q$ will contain also the belonging of control functions $q$ to $W^{3,\infty}(I)$ and the existence of a constant $c_\infty>0$ such that
	\beq
		\|q\|_{W^{3,\infty}}\leq c_\infty\qquad\forall\ q\in\Q.
	\label{eq:qinH3}
	\eeq
that is
\beq
	\Q := \{q\in W^{3,\infty}(I)\cap H^1_0(I)\colon q(x)\leq 1-\varepsilon,\ \forall x\in I, \text{ \ and \ } \|q\|_{W^{3,\infty}(I)}\leq c_\infty \}
\label{eq:QadW3}.
\eeq
Thanks to the above definition and to Lemma \ref{th:HkTq}, when handling with functions belonging to $H^k(\O)$ or $H^k(\Oq)$ for $k\leq3$, we can indifferently consider their norm in the physical domain $\Oq$ or in the reference domain $\O$.

Now we take into account the state problem, and we assume the validity of the following regularity result for Stokes problem:
\begin{ssmptn}
	Let $\Omega_q$ be an open, bounded set of $\mathds{R}^2$ and let $\Gamma_q$ be $C^{1,1}$ and $\partial\Oq\setminus\Gamma_q$ polygonal with $\partial\Oq$ having convex corners. Assume that data functions fulfill the following requests:
	\beq
		\nu\in H^3(\hat{\Omega}),\quad\eta\in H^2(\hat{\Omega}),\quad\mathbf f\in [H^2(\hat{\Omega})]^2,\quad\mathbf g_D\in[H^{7/2}(\Gamma_3)]^2,\quad\mathbf g_N\in[H^{5/2}(\Gamma_1)]^2,
	\eeq
	and suitable compatibility conditions.
	Then, for the solution $(\widetilde{\u},\widetilde{p})$ of \eqref{eq:Stokesdeb}, the following hold:
	\begin{enumerate}[\quad(a)]
		\item $(\widetilde{\u},\tilde{p})\in[H^3(\Oq)]^2\times H^2(\Oq)$
		\item $\|\nabla^3\tilde\u\|_\Oq+\|\nabla^2 \tilde p\|_\Oq \leq c(\eta,\nu,\mathbf g_D,\mathbf g_N,\mathbf f,\hat{\Omega}).$
	\end{enumerate}
\label{th:regolaritaS_app_lmm}
\end{ssmptn}
\begin{rmrk}
	Assumption \ref{th:regolaritaS_app_lmm} can be proved by resorting to results presented in \cite{Guo2006}. In that paper, weighted Sobolev spaces $H^{k,l}_\beta(\Omega)$ are considered, but the results can be brought back to classical Sobolev spaces exploiting the following inclusions: $H^k(\Omega_q)\subset H^{k,l}_\beta(\Omega_q)\subset H^{l-1}(\Omega_q)$, for any $k\geq l\geq 0$ and $\beta\in[0,1]^4$.
\end{rmrk}

The last ingredient that we need in order to prove a regularity result for the solution of our transformed problem \eqref{eq:StokesdebT} is represented by additional regularity requests on data. Since we want regularity not only for the solution of \eqref{eq:StokesdebT}, but also for its derivatives w.r.t.\ the control, namely $S'(q)(\dq),S''(q)(\dq,\dq)$, we have to assume a slightly stronger regularity of data than that considered in Assumption \ref{th:regolaritaS_app_lmm}.
\begin{ssmptn}
	Data functions have the following regularity:
	\beq
		\nu\in H^5(\hat{\Omega}),\quad\eta\in H^4(\hat{\Omega}),\quad\mathbf f\in [H^4(\hat{\Omega})]^2,\quad\mathbf g_D\in[H^{7/2}(\Gamma_3)]^2,\quad\mathbf g_N\in[H^{5/2}(\Gamma_1)]^2,
	\eeq
	and suitable compatibility conditions hold on data.
\label{ass:regolaritaS}
\end{ssmptn}

We are now ready to state a regularity result for the state variables and their shape derivatives.
\begin{thrm}
	Under Assumptions \ref{ass:regolaritaS}, \ref{th:regolaritaS_app_lmm}, there exist three positive constants $c_0,c_1,c_2$, such that for any $q\in\Q$, with $\|q\|_{W^{3,\infty}(I)}\leq c_\infty$, and for any $\dq,\tau q\in \delta Q$, and independently from them, it holds that
	\beq
		S(q),S'(q)(\dq),S''(q)(\dq,\tau q)\ \in\ [H^3(\Omega_0)]^2\times H^2(\Omega_0)\qquad\text{ and }
	\eeq
	\begin{enumerate}[\quad(a)]
		\item $\|S(q)\|_{[H^3(\Omega_0)]^2\times H^2(\Omega_0)}\leq c_0$
		\item $\|S'(q)(\dq)\|_{[H^3(\Omega_0)]^2\times H^2(\O)}\leq c_1\|\dq\|_{H^3(I)}$
		\item $\|S''(q)(\dq,\tau q)\|_{[H^3(\Omega_0)]^2\times H^2(\O)}\leq c_2\|\dq\|_{H^3(I)}\|\tau q\|_{H^3(I)}.$
	\end{enumerate}
\label{th:regolaritaS}
\end{thrm}
\begin{proof}
	Let $q\in\Q$, consider solution $(\u,p)=S(q)$ of the transformed problem \eqref{eq:StokesdebT} and remind that its physical counterpart $(\tilde\u,\tilde p)=\tilde S(q)$ is the solution of Stokes problem \eqref{eq:Stokesdeb} on $\Oq$.\\
	Now, we can verify the hypotheses of Lemma \ref{th:regolaritaS_app_lmm}: $\Oq$ is surely an open bounded subset of $\mathds{R}^2$; its boundary $\Gamma_q$ is $C^{1,1}$ because it is the graph of the control function $q\in\Q\subset H^3(I)\subset C^{1,1}(\overline{\Oq})$ and for the same reason its terminal points cannot present a concave angle; the regularity of external force and boundary data, together with the compatibility conditions, are given by Assumption \ref{ass:regolaritaS}. Then, Lemma \ref{th:regolaritaS_app_lmm} holds and we have $(\tilde\u,\tilde p)\in [H^3(\Oq)]^2\times H^2(\Oq)$ and $\|\nabla^3\tilde{\u}\|_\Oq+\|\nabla^2 \tilde{p}\|_\Oq \leq c(\mathbf f,\mathbf g_D,\mathbf g_N, \hat{\Omega})$.
Finally, the results on $(\tilde\u,\tilde p)$ directly transfer to $(\u,p)$, thanks to Lemma \ref{th:HkTq}.

For points (\emph{b}) and (\emph{c}), the proof is exactly the same, considering Assumption \ref{ass:regolaritaS} in order to control the more complex right-hand sides appearing dealing with $S'(q)(\dq)$ and $S''(q)(\dq,\dq)$, with the aim of prove the validity of the hypotheses of Lemma \ref{th:regolaritaS_app_lmm}. The dependence on $\|\dq\|_{H^3(I)}$ on the right-hand side comes out from the bounds of the coefficients, similar to those reported in Proposition \ref{th:coeff}.

\end{proof}

So far we have obtained a regularity result for the state variables: now we want to show that the optimal control belongs to $H^5(I)$. Indeed, this regularity holds for any $q\in\Q$ satisfying the first order optimality condition, as stated in the following result:
\begin{thrm}
	Let $\overline{q}\in\Q$ be such that optimality condition \eqref{eq:ottimalita} holds in $\overline{q}$.
	Then it holds that $\overline{q}\in H^5(I)$.
\label{th:regolaritaq}
\end{thrm}
\begin{proof}
Let us take into account Hadamard formula for $j'$, given by Lemma \ref{th:Hadamard}, i.e. $j'(q)(\dq)= 2\alpha(q'',\dq'')_I + (\Psi,\dq)_I$.
We start by noticing that $\beta(\int_Iq(x)dx - \overline{V})$ is constant, then certainly belonging to $H^1(I)$.\\
The regularity Theorem \ref{th:regolaritaS} can be applied to both the state variables $(\u,p)$ and the adjoint state variables $(\z,s)$. Then, thanks to the definition \eqref{eq:qinH3} and Lemma \ref{th:HkTq}, we get
\beq
	(\widetilde{\u},\widetilde{\z},\widetilde{s}) = (\u,\z,s)\circ T_q\quad \in [H^3(\Oq)]^2\times[H^3(\Oq)]^2\times H^2(\Oq).
\eeq
Taking the traces of $\nabla{\widetilde{\u}},\nabla{\widetilde{\z}},\widetilde{s}$ on boundary $\Gamma_q$ and using its parametrization $\gamma:x\mapsto(x,q(x))$, we get
\beq
	(\nabla\widetilde{\u},\nabla\widetilde{\z},\widetilde{s})(x,q(x)) \quad \in \quad [H^{3/2}(I)]^{2\times2}\times[H^{3/2}(I)]^{2\times2}\times H^{3/2}(I).
\label{eq:uzs_psireg}
\eeq
Thanks to this regularity, together with the continuous embedding $H^{3/2}(I)\hookrightarrow W^{1,4}(I)$, we can conclude that $\Psi$ belongs to $H^1(I)$.

Now, taking a $\overline{q}\in\Q$ such that the optimality condition $j'(\overline{q})(\dq)=0$ holds, we get
\beq
	\int_I\overline{q}''\dq'' dx = -\int_I \frac{1}{2\alpha}\Psi \dq\, dx \qquad \forall \dq\in C^\infty_0(\overline{I}).
\label{eq:der4}
\eeq
Finally, we observe that \eqref{eq:der4} is equivalent to say that the fourth weak derivative of $\overline{q}$ is exactly $-\frac{1}{2\alpha}\Psi$. Being $\alpha$ a non-zero constant and belonging $\Psi$ to $H^1(I)$, we get $\overline{q}^{(iv)}\in H^1(I)$. Since we already have $\overline{q}\in W^{3,\infty}\subset H^3(I)$ (see \eqref{eq:qinH3}), we obtain the thesis, i.e. $\overline{q}\in H^5(I)$.
\end{proof}

\section{Results for the coercivity of functional \texorpdfstring{$j$}{j}}\label{sec:per314}

In this appendix, we present two useful results for the proof of Lemma \ref{th:314}. The first concerns the sequential continuity of the state operator derivatives w.r.t.\ the variations of control.
\begin{lmm}
	Let $q\in\Q$ and consider a sequence $\{\dq_n\}_{n\in\mathds{N}}\subset Q$. If there exists a $\dq\in Q$ such that $\dq_n\to\dq$ in $C^1(\overline{I})$, then
	\begin{enumerate}[(a)]
	\item $S'(q)(\dq_n)\to S'(q)(\dq)$ in $V\times P$
	\item $S''(q)(\dq_n,\dq_n)\to S''(q)(\dq,\dq)$ in $V\times P$
	\end{enumerate}
\label{th:KinVex39}
\end{lmm}
\begin{proof}
	Because of the linearity and the well-posedness of problems \eqref{eq:dS}, \eqref{eq:d2S}, it suffices to prove the convergence of the right-hand sides in $V'\times P'$: this is obtained from the continuity of $\dot{F},\dot{a},\dot{b},$ $\ddot{F},\ddot{a},\ddot{b}$ w.r.t.\ the variation $\dq$.\\
	We just give an example of the steps to be taken, processing a term from $\ddot{a}(q,\dq\dq)(\u,\v)$:
	\begin{equation*}\begin{split}
	&\left|\left(\nabla\nu^q\cdot V_{\dq_n}\nabla\u\,A'_{q,\dq_n},\nabla\v\right)-\left(\nabla\nu^q\cdot V_\dq\nabla\u\,A'_{q,\dq},\nabla\v\right)\right|\leq\\
	&\leq\left|\left(\nabla\nu^q\cdot V_{\dq_n}\nabla\u(A'_{q,\dq_n}-A'_{q,\dq}),\nabla\v\right)\right|+\left|\left(\nabla\nu^q\cdot (V_{\dq_n}-V_\dq)\nabla\u\,A'_{q,\dq},\nabla\v\right)\right|\leq\\
	&\leq\|\nu\|_{W^{1,\infty}(\hat{\Omega})}\|\nabla\u\|\|\nabla\v\|\left(\|A'_{q,\dq_n}-A'_{q,\dq}\|_\infty\|V_{\dq_n}\|_\infty+\|V_{\dq_n}-V_\dq\|_\infty\|A'_{q,\dq}\|_\infty\right)
	\end{split}\end{equation*}
	The convergence of $\dq_n$ in $C^1(\overline{I})$ implies the uniform convergence of $A'_{q,\dq_n}-A'_{q,\dq}$ and $V_{\dq_n}-V_\dq$ to zero, as it can be seen from the definition of such quantities. Moreover, being $\{\dq_n\}$ bounded in $C^1(\overline{I})$, Proposition \ref{th:coeff} ensures that $\|V_{\dq_n}\|_\infty, \|A'_{q,\dq}\|_\infty$ are bounded themselves.
\end{proof}
From Lemma \ref{th:KinVex39}, using the Dominated Convergence Theorem and the compact embedding $H^2(I)\subset\subset C^1(\overline{I})$ yields a similar result for the derivatives of cost functional $j$:
\begin{crllr}
	Let $q\in\Q$ and $\{\dq_n\}_{n\in\mathds{N}}\subset\delta Q$ such that there exists a $\dq\in\delta Q$ for which $\dq_n\weak\dq$ in $H^2(I)$. Then,
	\beq
		j'(q)(\dq_n)\underset{n\to\infty}{\longrightarrow} j(q)(\dq),\qquad j''(q)(\dq,\dq)\leq \liminf_{n\to\infty}j''(q)(\dq_n,\dq_n).
	\eeq
\label{th:KinVex310}
\end{crllr}

\end{appendices}
\bibliography{library}

\end{document}